\documentclass[preprint,12pt,draft]{elsarticle}
\usepackage{amsthm,amssymb,amsfonts,amsmath,times1}
\usepackage{color}
\usepackage{hyperref}
\usepackage{mathrsfs}

\numberwithin{equation}{section}

\theoremstyle{plain}
\newtheorem{theorem}{Theorem}[section]
\newtheorem{proposition}[theorem]{Proposition}
\newtheorem{lemma}[theorem]{Lemma}
\newtheorem{corollary}[theorem]{Corollary}

\theoremstyle{remark}
\newtheorem{remark}{Remark}[section]
\newtheorem{example}{Example}[section]

\theoremstyle{definition}
\newtheorem{definition}{Definition}[section]

\newcommand{\Div}{\mathop{\mathrm{div}}\nolimits}
\newcommand{\ext}{\mathop{\mathrm{ext}}\nolimits}
\newcommand{\Ker}{\mathop{\mathrm{Ker}}\nolimits}
\newcommand{\supp}{\mathop{\mathrm{supp}}\nolimits}
\newcommand{\Diff}{\mathrm{Diff}}
\newcommand{\Vect}{\mathrm{Vect}}
\newcommand{\const}{\mathrm{const}}
\newcommand{\rd}{\mathrm{d}}
\newcommand{\re}{\mathrm{e\myp}}

\newcommand{\CD}{\mathop{\raisebox{-.11pc}{\Large$\cdot$}}}

\newcommand{\RR}{\mathbb{R}}
\newcommand{\NN}{\mathbb{N}}
\newcommand{\ZZ}{\mathbb{Z}}

\newcommand{\calB}{\mathcal{B}}

\newcommand{\frakX}{\mathfrak{X}}

\newcommand{\g}{\text{\slshape\textrm{g}\myp
}}
\newcommand{\hatg}{\hat{\text{\slshape\textrm{g}\mypp}}\myn{}}
\newcommand{\gcl}{\g_{\myn\mathrm{cl}}}

\newcommand{\hatgamma}{{\hat{\gamma\myp}\myn{}}}
\newcommand{\Q}{\mathcal{I}_{\mathfrak{q}}}
\newcommand{\U}{{\varPhi}}
\newcommand{\GR}{{\mathscr{G}_{\mathrm{R}}}}

\newcommand{\myp}{\mbox{$\:\!$}}
\newcommand{\mypp}{\mbox{$\;\!$}}

\newcommand{\myn}{\mbox{$\;\!\!$}}
\newcommand{\mynn}{\mbox{$\:\!\!$}}

\begin{document}

\begin{frontmatter}

\title{Gibbs cluster measures on configuration spaces}


\author[Leeds]{Leonid Bogachev\fnref{t1,t2}\corref{cor1}}
\ead{L.V.\,Bogachev@leeds.ac.uk}
\author[York]{Alexei Daletskii\fnref{t2}}
\ead{ad557@york.ac.uk}

\address[Leeds]{Department of Statistics, University of Leeds, Leeds LS2 9JT, UK}
\address[York]{Department of Mathematics, University of York, York YO10 5DD, UK}


\fntext[t1]{Research supported in part by a Leverhulme Research
Fellowship.}
\fntext[t2]{Research supported in part by DFG Grant
436\,RUS\,113/722.}
\cortext[cor1]{Corresponding author.}

\begin{abstract}
The distribution $\gcl$ of a Gibbs cluster point process in
$X=\mathbb{R}^{d}$ (with i.i.d.\ random clusters attached to points
of a Gibbs configuration with distribution $\g$) is studied via the
projection of an auxiliary Gibbs measure $\hatg$ in the space of
configurations $\hatgamma=\{(x,\bar{y})\}\subset
X\times\mathfrak{X}$, where $x\in X$ indicates a cluster ``center''
and $\bar{y}\in\mathfrak{X}:=\bigsqcup_{\myp n}\myn X^n$ represents
a corresponding cluster relative to $x$.
We show that the measure $\gcl$ is quasi-invariant with respect to
the group $\Diff_{0}(X)$ of compactly supported diffeomorphisms of
$X$, and prove an integration-by-parts formula for $\gcl$. The
associated equilibrium stochastic dynamics is then constructed using
the method of Dirichlet forms. These results are quite general; in
particular, the uniqueness of the background Gibbs measure $\g$ is
not required. The paper is an extension of the earlier results for
Poisson cluster measures
[J.~Funct.\ Analysis 256 (2009) 432--478], where a different
projection construction was utilized
specific to this ``exactly soluble'' case.
\end{abstract}


\end{frontmatter}

\section{Introduction}\label{sec:1}

The concept of particle configurations is instrumental in
mathematical modelling of multi-component stochastic systems. Rooted
in statistical mechanics and theory of point processes, the
development of the general mathematical framework for suitable
classes of configurations has been
a recurrent research theme fostered by widespread applications
%
across the board, including
quantum physics, astrophysics, chemical physics, biology, ecology,
computer science, economics, finance, etc.\ (see an extensive
bibliography in \cite{DVJ1}).

In the past 15 years or so, there has been a more specific interest
in the \emph{analysis} on configuration spaces. To fix basic
notation, let $X$ be a topological space (e.g., a Euclidean space
$X=\RR^{d}$), and let $\varGamma_X=\{\gamma\}$ be the configuration
space over $X$, that is, the space of countable subsets (called
\emph{configurations}) $\gamma\subset X$ without accumulation
points. Albeverio, Kondratiev and R\"{o}ckner \cite{AKR1,AKR2} have
proposed an approach to configuration spaces $\varGamma_X$ as
\emph{infinite-dimensional manifolds}, based on the choice of a
suitable probability measure $\mu $ on $\varGamma_{X}$
which is
quasi-invariant with respect to $\Diff_{0}(X)$, the group of
compactly supported diffeomorphisms of $X$. Providing that the
measure $\mu$ can be shown to satisfy an integration-by-parts
formula, one can construct, using the theory of Dirichlet forms, an
associated equilibrium dynamics (stochastic process) on
$\varGamma_{X}$ such that $\mu $ is its invariant measure
\cite{AKR1,AKR2,MR} (see
\cite{ADKal,AKR2,AKR3,
Ro} and references therein for further discussion of various
theoretical aspects and applications).

This general programme has been first implemented in \cite{AKR1} for
the \textit{Poisson} measure $\mu$ on $\varGamma_{X}$, and then
extended in \cite{AKR2} to a wider class of \textit{Gibbs} measures,
which appear in statistical mechanics of classical continuous gases.
In the Poisson case, the canonical equilibrium dynamics is given by
the well-known independent particle process, that is, an infinite
family of independent (distorted) Brownian motions started at the
points of a random Poisson configuration. In the Gibbsian case, the
equilibrium dynamics is much more complex due to interaction between
the particles.

In our earlier papers \cite{BD1,BD3}, a similar analysis was
developed for a different class of random spatial structures, namely
\emph{Poisson cluster point processes}, featured by spatial grouping
(``clustering'') of points around the background random (Poisson)
configuration of invisible ``centers''. Cluster models are well
known in the general theory of random point processes \cite{CI,DVJ1}
and are widely used in numerous applications ranging from
neurophysiology (nerve impulses) and ecology (spatial aggregation of
species) to seismology (earthquakes) and cosmology (constellations
and galaxies); see \cite{BD3,CI,DVJ1} for some references to
original papers.


Our technique in \cite{BD1,BD3} was based on the representation of a
given Poisson cluster measure on the configuration space
$\varGamma_{X}$ as the projection image of an auxiliary Poisson
measure on a more complex configuration space
$\varGamma_{\mathfrak{X}}$ over the disjoint-union space
$\mathfrak{X}:=\bigsqcup_{\myp n} \mynn X^{n}$,
with ``droplet'' points $\bar{y}\in \mathfrak{X}$ representing
individual clusters (of variable size).
%
The principal advantage of this construction is that it allows one
to apply the well-developed apparatus of Poisson measures to the
study of the Poisson cluster measure.

In the present paper,\footnote{Some of our results have been
announced in
\cite{BD2} (in the case of clusters of
fixed size).} our aim is to extend this approach to a more general
class of \emph{Gibbs cluster measures} on the configuration space
$\varGamma_X$, where
the distribution of cluster centers is given by a Gibbs (grand
canonical) measure $\g\in\mathscr{G}(\theta,\U\myp)$ on
$\varGamma_X$, with a reference measure $\theta$ on $X$ and an
interaction potential $\U$.
We focus on Gibbs cluster processes in $X=\mathbb{R}^{d} $ with
i.i.d.\ random clusters of random size.
Let us point out that we do not require the uniqueness of the Gibbs
measure, so our results are not affected by possible ``phase
transitions'' (i.e., non-uniqueness of
$\g\in\mathscr{G}(\theta,\U\myp)$). Under some natural smoothness
conditions on the reference measure $\theta$ and the distribution
$\eta$ of the generic cluster, we prove the
$\Diff_{0}(X)$-quasi-invariance of the corresponding Gibbs cluster
measure $\gcl$ (Section~\ref{sec:3.2}), establish the
integration-by-parts formula (Section~\ref{sec:3.3}) and construct
the associated Dirichlet operator, which leads to the existence of
the equilibrium stochastic dynamics on the configuration space
$\varGamma_{X}$ (Section~\ref{sec:4}).

Unlike the Poisson cluster case, it is now impossible to work with
the measure arising in the space $\varGamma_{\mathfrak{X}}$ of
droplet configurations $\bar{\gamma}=\{\bar{y}\}$, which is hard to
characterize for Gibbs cluster measures. Instead, in order to be
able to pursue our projection approach while still having a
tractable pre-projection measure, we choose the configuration space
$\varGamma_{\mathcal{Z}}$ over the set
$\mathcal{Z}:=X\times\mathfrak{X}$, where each configuration
$\hatgamma\in \varGamma_{\mathcal{Z}}$ is a (countable) set of pairs
$z=(x,\bar{y})$ with $x\in X$ indicating a cluster center and
$\bar{y}\in \mathfrak{X}$ representing a cluster attached to $x$. A
crucial step is to show that the corresponding measure $\hatg$ on
$\varGamma_{\mathcal{Z}}$ is again Gibbsian, with the
reference measure $\sigma =\theta\otimes \eta $ and a ``cylinder''
interaction potential \,$\hat{\U}(\hatgamma):=
\U(\mathfrak{p}(\hatgamma))$, where $\varPhi$ is the original
interaction potential associated with the background Gibbs measure
$\g$ and $\mathfrak{p}$ is the operator on the configuration space
$\varGamma_{\mathcal{Z}}$ projecting a configuration
$\hatgamma=\{(x,\bar{y})\}$ to the configuration of cluster centers,
$\gamma=\{x\}$. We then project the Gibbs measure $\hatg$ from the
``higher floor'' $\varGamma_{\mathcal{Z}}$ directly to the
configuration space $\varGamma_{X}$ (thus skipping the
``intermediate floor'' $\varGamma_{\mathfrak{X}}$), and show that
the resulting measure coincides with the original Gibbs cluster
measure $\gcl$ (Section~\ref{sec:2}).

In fact, it can be we shown (Section \ref{sec:2.3}) that
\textit{any} cluster measure $\mu_{\rm cl}$ on $\varGamma_X$ can be
obtained by a similar projection from $\varGamma_{\mathcal{Z}}$.
Even though it may not always be possible to find an intrinsic
characterization of the corresponding lifted measure $\hat{\mu}$ on
the configuration space $\varGamma_{\mathcal{Z}}$ (unlike the
Poisson and Gibbs cases), we expect that the projection approach can
be instrumental in the study of more general cluster point processes
by a reduction to point processes in more complex phase spaces but
with a simpler correlation structure. We intend to develop these
ideas elsewhere.

\section{Gibbs cluster measures via projections} 
\label{sec:2}

In this section, we start by recalling some basic concepts and
notations for random point processes and associated probability
measures in configuration spaces (Section~\ref{sec:2.1}), followed
in Section~\ref{sec:2.2} by a definition of a general cluster point
process (CPP). In Section~\ref{sec:2.3}, we explain our main
``projection'' construction allowing one to represent CPPs in the
phase space $X$ in terms of auxiliary measures on a more complex
configuration space involving Cartesian powers of $X$. The
implications of such a description are discussed in greater detail
for the particular case of Gibbs CPPs (Sections
\ref{sec:2.4},~\ref{sec:2.5}).

\subsection{Probability measures on configuration spaces}\label{sec:2.1}

Let $X$ be a Polish space
equipped with the Borel $\sigma$-algebra ${\mathcal{B}}(X)$
generated by the open sets. Denote $\overline{\ZZ}_+:=
\ZZ_+\cup\{\infty\}$, where $\ZZ_{+}:=\{0,1,2,\dots \}$, and
consider a space ${\mathfrak{X}}$ built from all Cartesian powers of
$X$, that is, the disjoint union
\begin{equation}\label{eq:calX}
\frakX:={\textstyle\bigsqcup\limits_{\myp n\in\overline{\ZZ}_{+}}}
X^{n},
\end{equation}
including $X^{0}=\{\emptyset \}$ and the space $X^\infty$ of
infinite sequences $(x_1, x_2,\dots)$. That is,
$\bar{x}=(x_{1},x_{2},\dots )\in{\mathfrak{X}}$ if and only if
$\bar{x}\in X^{n}$ for some $n\in\overline{\ZZ}_{+}$.
We take the liberty to write $x_{i}\in\bar{x}$ if
$x_{i}$ is a coordinate of the ``vector'' $\bar{x}$.
%
The space $\mathfrak{X}$ is endowed with the natural disjoint union
topology induced by the topology in $X$.

\begin{remark}\label{rm:compact}
Note that a set $K\subset{\mathfrak{X}}$ is compact if and only if
$K=\bigsqcup_{\myp n=0}^{N} K_{n}$, where $N<\infty $ and $K_{n}$
are compact subsets of $X^{n}$, respectively.
\end{remark}

\begin{remark}
$\mathfrak{X}$ is a Polish space as a disjoint union of Polish spaces.
\end{remark}

Denote by ${\mathcal{N}}(X)$ the space of
$\overline{\ZZ}_{+}$-valued measures $N(\cdot)$ on
${\mathcal{B}}(X)$ with countable (i.e., finite or countably
infinite) support.
Consider the natural projection
\begin{equation}\label{eq:pr0}
{\mathfrak{X}}\ni \bar{x}\mapsto
{\mathfrak{p}}(\bar{x}):=\sum_{x_{i}\in \bar{x}}\delta_{x_{i}}\in
{\mathcal{N}}(X),
\end{equation}
where $\delta_{x}$ is the Dirac measure at point $x\in X$.
That is to say, under the map ${\mathfrak{p}}$ each vector from
${\mathfrak{X}}$ is ``unpacked'' into its components to yield a
countable aggregate of (possibly multiple) points in $X$, which can
be interpreted as a generalized configuration $\gamma$,
%
\begin{equation}\label{eq:pr}
\mathfrak{p}(\bar{x})\leftrightarrow\gamma
:={\textstyle\bigsqcup\limits_{x_{i}\in \bar{x}}}\{x_{i}\},\qquad
\bar{x}=(x_{1},x_{2},\dots )\in {\mathfrak{X}}.
\end{equation}

In what follows, we
interpret the notation $\gamma $ either as an aggregate of
points in $X$ or as a $\overline{\mathbb{Z}}_{+}$-valued measure or
both, depending on the context. Even though generalized
configurations are not, strictly speaking, subsets of $X$ (because
of possible multiplicities), it is convenient to
use set-theoretic notations, which should not cause any confusion.
For instance, we write $\gamma \cap B$ for the restriction of
configuration $\gamma $ to a subset $B\in {\mathcal{B}}(X)$.
For a function $f\myn:X\to\RR$ we denote
\begin{equation}\label{eq:f-gamma}
\langle f,\gamma \rangle :=\sum_{x_i\in \gamma}f(x_i)
\equiv\int_{X}f(x)\,\gamma({\rd}x).
\end{equation}
In particular, if $\mathbf{1}_{B}(x) $ is the indicator function of
a set $B\in{\mathcal{B}}(X)$ then $\langle \mathbf{1}_{B},\gamma
\rangle =\gamma (B)$ is the total number of points (counted with
their multiplicities) in $\gamma\cap B$.

\begin{definition}\label{def:gen}
A \textit{configuration space} $\varGamma_{X}^{\sharp}$ is the set
of generalized configurations $\gamma$ in $X$, endowed with the
\textit{cylinder $\sigma$-algebra}
$\mathcal{B}(\varGamma_{X}^{\sharp})$ generated by the class of
cylinder sets $C_{B}^{\myp n}:=\{\gamma\in\varGamma_{X}^{\sharp}:
\gamma (B)=n\}$, \,$B\in{\mathcal{B}}(X)$, \,$n\in\ZZ_+$\myp.
\end{definition}
\begin{remark}
It is easy to see that the map
$\mathfrak{p}:{\mathfrak{X}}\rightarrow \varGamma_{X}^{\sharp} $
defined by formula (\ref{eq:pr}) is measurable.
\end{remark}

In fact, conventional theory of point processes (and their
distributions as probability measures on configuration spaces)
usually rules out the possibility of accumulation points or multiple
points (see, e.g., \cite{DVJ1}).

\begin{definition}\label{def:proper}
A configuration $\gamma \in \varGamma_{X}^{\sharp}$ is said to be
\emph{locally finite} if $\gamma(B)<\infty $ for any compact set
$B\subset X$. A configuration $\gamma \in \varGamma_{X}^{\sharp}$ is
called \textit{simple} if $\gamma(\{x\})\le 1$ for each $x\in X$. A
configuration $\gamma \in \varGamma_{X}^{\sharp}$ is called
\emph{proper} if it is both locally finite and simple. The set of
proper configurations will be denoted by $\varGamma_{X}$ and called
the \textit{proper configuration space} over $X$. The corresponding
$\sigma$-algebra ${\mathcal{B}}(\varGamma_{X})$ is generated by the
cylinder sets $\{\gamma \in \varGamma_{X}:\gamma (B)=n\}$ \,($B\in
{\mathcal{B}}(X)$, \,$n\in\ZZ_+$).
\end{definition}

%
Like in the standard theory based on proper configuration spaces
(see, e.g., \cite[\S\,6.1]{DVJ1}), every probability measure $\mu$
on the generalized configuration space $\varGamma_{X}^{\sharp}$ can
be characterized by its Laplace functional (cf.\ \cite{BD3})
\begin{equation}
L_{\mu}(f):=
\int_{\varGamma_{X}^\sharp} \re^{-\langle
f,\myp\gamma \rangle}\,\mu(\rd \gamma),\qquad f\in
{\mathrm{M}}_{+}(X), \label{eq:LAPLACE}
\end{equation}
where ${\mathrm{M}}_{+}(X)$ is the class of measurable non-negative
functions on $X$.

\subsection{Cluster point processes}
\label{sec:2.2}

Let us recall the notion of a general cluster point process (CPP).
Its realizations are constructed in two steps: (i) a background
random configuration of (invisible) ``centers'' is obtained as a
realization of some point process $\gamma_{\mathrm{c}}$ governed by
a probability measure $\mu_{\mathrm{c}}$ on $\varGamma_{X}^{\sharp
}$, and (ii) relative to each center $x\in \gamma_{\mathrm{c}}$, a
set of observable secondary points (referred to as a \emph{cluster}
centered at~$x$) is generated according to a point process
$\gamma_{x}^{\myp\prime}$ with probability measure $\mu_{x}$ on
$\varGamma_{X}^{\sharp }$ ($x\in X$).

The resulting (countable) assembly of random points, called the
\emph{cluster point process}, can be symbolically expressed as
\begin{equation*}
\gamma ={\textstyle\bigsqcup\limits_{x\in \gamma_{\mathrm{c}}}}\myp
\gamma_{x}^{\myp\prime }\in \varGamma_{X}^{\sharp},
\end{equation*}
where the disjoint union signifies that multiplicities of points
should be taken into account. More precisely, assuming that the
family of secondary processes $\gamma_x^{\myp\prime}(\cdot)$ is
measurable as a function of $x\in X$, the integer-valued measure
corresponding to a CPP realization $\gamma $ is given by
\begin{equation}  \label{eq:cluster-gamma}
\gamma(B)=\int_{X}\gamma_x^{\myp\prime}(B)\,\gamma_{\mathrm{c}}(\rd
x) =\sum_{x\in
\gamma_{\mathrm{c}}}\gamma_x^{\myp\prime}(B),
\qquad B\in {\mathcal{B}}(X).
\end{equation}

In what follows, we assume
that (i) $X$ is a linear space (e.g., $X=\RR^d$) so that
translations $X\ni y\mapsto y+x\in X$ are defined, and (ii) random
clusters are independent and identically distributed (i.i.d.), being
governed by the same probability law translated to the cluster
centers, so that, for any $x\in X$, we have
$\mu_{x}(A)=\mu_{0}(A-x)$ ($A\in\mathcal{B}(\varGamma_X^{\sharp})$.

%
In turn, the measure $\mu_{0}$ on $\varGamma_X^{\sharp}$
determines a probability distribution $\eta $ in ${\mathfrak{X}}$
which is symmetric with respect to permutations of coordinates.
Conversely, $\mu_{0}$ is a push-forward of the measure $\eta $ under
the projection map ${\mathfrak{p}}:{\mathfrak{X}}\rightarrow
\varGamma_{X}^{\sharp }$ defined by (\ref{eq:pr}), that is,
\begin{equation}\label{eq:p*eta}
\mu_{0}={\mathfrak{p}}^{*}\eta \equiv \eta \circ
{\mathfrak{p}}^{-1}.
\end{equation}

\begin{remark}
Unlike the standard CPP theory when sample configurations are
\emph{presumed} to be a.s.\ locally finite (see, e.g.,
\cite[Definition 6.3.I]{DVJ1}), the description of the CPP given
above only implies that its configurations $\gamma $ are countable
aggregates in $X$, but possibly with multiple and/or accumulation
points, even if the background point process $\gamma_{\mathrm{c}}$
is proper. Therefore, the distribution $\mu $ of the CPP
(\ref{eq:cluster-gamma}) is a probability measure defined on the
space $\varGamma_{X}^{\sharp }$ of \emph{generalized}
configurations. It is a matter of interest to obtain conditions in
order that $\mu $ be actually supported on the proper configuration
space $\varGamma_{X}$, and we shall address this issue in Section
\ref{sec:2.4} below for Gibbs CPPs (see \cite{BD3} for the case of
Poisson CPPs).
\end{remark}

%

The following fact is well known in the case of CPPs without
accumulation points (see, e.g., \cite[\S \,6.3]{DVJ1}); its proof in
the general case is essentially the same
(see \cite[Proposition 2.5]{BD3}).

\begin{proposition}\label{pr:cluster}
Let $\mu_{{\rm cl}}$ be a probability measure on
$(\varGamma_X^{\sharp},\mathcal{B}(\varGamma_X^\sharp))$ determined
by the probability distribution of a CPP
\textup{(\ref{eq:cluster-gamma})}. Then its Laplace functional
is given, for all functions $f\in {\mathrm{M}}_+(X)$, by
\begin{equation}\label{laplace-G}
L_{\mu_{{\rm cl}}}(f)
=\int_{\varGamma^\sharp_{X}}\prod_{x\in \gamma_{\rm c}}\biggl(
\int_{\mathfrak{X} }\exp \biggl(- \sum_{y_i\in
\bar{y}}f(y_i+x)\Bigr) \,\eta ({\rd}\bar{y})\biggr) \,\mu_{\rm
c}({\rd}\gamma_{\rm c}).
\end{equation}
\end{proposition}

\subsection{A projection construction of cluster measures on configurations}\label{sec:2.3}

Denote $\mathcal{Z}:=X\times\frakX$,
and consider the space
$\varGamma_{\mathcal{Z}}^\sharp=\{\hatgamma\}$ of (generalized)
configurations
in $\mathcal{Z}$.
Let $p_{X}\mynn:\mathcal{Z}\to X$ be the natural projection to the
first coordinate,
\begin{equation}\label{eq:pX0}
\mathcal{Z}\ni z=(x,\bar{y})\mapsto p_{X}(z):=x\in X,
\end{equation}
and consider its pointwise lifting to the configuration space
$\varGamma^\sharp_{\mathcal{Z}}$ (preserving the same notation
$p_{X}$), defined as follows
\begin{equation}\label{eq:pX}
\varGamma_{\mathcal{Z}}^\sharp\ni \hatgamma\mapsto
p_{X}(\hatgamma):={\textstyle\bigsqcup\limits_{z\in\myp\hat{\gamma\myp}}}\{p_{X}(z)\}
\in\varGamma^\sharp_X.
\end{equation}

Let $\mu_{\rm cl}$ denote the probability measure on the
configuration space $\varGamma^\sharp_{X}$ associated with an
i.i.d.\ cluster point process (see Section~\ref{sec:2.2}), specified
by measures $\mu_{\rm c}$ on $\varGamma^\sharp_{X}$ and
$\eta$ on $\frakX$.

\begin{definition}\label{def:muhat}
Let us define a probability measure $\hat{\mu}$ on
$\varGamma_{\mathcal{Z}}^\sharp$ as the distribution of random
configurations $\hatgamma$ over $\mathcal{Z}$ obtained from
configurations $\gamma\in\varGamma^\sharp_{X}$ by attaching to each
point $x\in\gamma$ an i.i.d.\ random vector $\bar{y}_x$ with
distribution $\eta$:
\begin{equation}\label{eq:gamma-hat}
\varGamma^\sharp_{X}\ni\gamma_{\rm c}
\mapsto \hatgamma:={\textstyle\bigsqcup\limits_{x\in\gamma_{\rm
c}}\{(x,\bar{y}_x)\}} \in \varGamma_{\mathcal{Z}}^\sharp.
\end{equation}
Geometrically, the construction (\ref{eq:gamma-hat}) may be viewed
as random i.i.d.\ pointwise translations of configurations $\gamma$
from $X$ into the ``plane'' $\mathcal{Z}=X\times \mathfrak{X}$. The
measure $\hat{\mu}$ so obtained may be expressed in the differential
form as a skew product
\begin{equation}\label{eq:muhat0}
\hat{\mu}(\rd\hatgamma)=\mu_{\rm
c}(p_{X}(\rd\hatgamma))\,{\textstyle\bigotimes\limits_{z\in\hat{\gamma\myp}}}
\,\eta(p_\frakX(\rd{z})), \qquad
\hatgamma\in\varGamma^\sharp_{\mathcal{Z}}.
\end{equation}
Equivalently, for any function
$F\in{\mathrm{M}}_+(\varGamma^\sharp_{\mathcal{Z}})$,
\begin{equation}\label{eq:muhat}
\int_{\varGamma_{\mathcal{Z}}^\sharp}
F(\hatgamma)\,\hat{\mu}(\rd\hatgamma)
=\int_{\varGamma_{X}^\sharp}\!\biggl(\int_{{\mathfrak{X}}^\infty}
F\Bigl(\mypp{\textstyle\bigcup\limits _{x\in\gamma}}\{({\textstyle
x},\bar{y})\}\Bigr)\,{\textstyle\bigotimes\limits_{x\in\gamma}}
\,\eta(\rd\bar{y})\mynn\biggr)\,\mu_{\rm c}(\rd\gamma).
\end{equation}
\end{definition}

\begin{remark}
Note that formula
(\ref{eq:muhat}) is a simple case of the general disintegration
theorem, or the ``total expectation formula''
(see, e.g., \cite[Theorem 5.4]{Kal}
or
\cite[Ch.\,V, \S\mypp8, Theorem~8.1]{Par}).
\end{remark}



Recall that the ``unpacking'' map
$\mathfrak{p}:\frakX\to\varGamma_X^\sharp$ is defined in
(\ref{eq:pr}), and consider a map
$\mathfrak{q}:\mathcal{Z}\rightarrow \varGamma_{X}^{\sharp}$ acting
by the formula
\begin{equation}\label{proj1}
\mathfrak{q}(x,\bar{y}):=\mathfrak{p}(\bar{y}+x)=
{\textstyle\bigsqcup\limits_{y_i\in\bar{y}}}\{y_i+x\},\qquad
(x,\bar{y})\in \mathcal{Z}.
\end{equation}
Here and below, we use the shift notation ($x\in X$)
\begin{equation}\label{eq:shift}
\bar{y}+x:=(y_1+x,\,y_2+x,\dots),\qquad \bar{y}=(y_1\myn,
y_2,\dots)\in \frakX,
\end{equation}
In the usual ``diagonal'' way, the map $\mathfrak{q}$ can be lifted
to the configuration space $\varGamma_{\mathcal{Z}}^{\sharp}$:
\begin{equation}\label{eq:proj}
\varGamma_{\mathcal{Z}}^{\sharp}\ni \hatgamma\mapsto
\mathfrak{q}(\hatgamma):={\textstyle\bigsqcup\limits_{z\in\myp
\hat{\gamma\myp}}}\mathfrak{q}(z)\in \varGamma_{X}^\sharp.
\end{equation}

\begin{proposition}\label{pr:q-meas}
The map $\mathfrak{q}:\varGamma_{\mathcal{Z}}^{\sharp} \to
\varGamma_{X}^{\sharp}$ defined by \textup{(\ref{eq:proj})} is
measurable.
\end{proposition}
\proof Observe that $\mathfrak{q}$ can be represented as a
composition
\begin{equation}\label{eq:circ}
\mathfrak{q}=\mathfrak{p}\circ p_{\mathfrak{X}}\circ
p_+:\varGamma_{\mathcal{Z}}^\sharp\stackrel{p_+}{\longrightarrow}
\varGamma_{\mathcal{Z}}^\sharp\stackrel{p_{\mathfrak{X}}}{\longrightarrow}
\varGamma_{\mathfrak{X}}^\sharp\stackrel{\mathfrak{p}}{\longrightarrow}
\varGamma_{X}^\sharp,
\end{equation}
where the maps $p_+$\myp, $p_{\mathfrak{X}}$ and $\mathfrak{p}$ are
defined, respectively, by
\begin{align}
\label{eq:p+} \varGamma_{\mathcal{Z}}^{\sharp}\ni \hatgamma&\mapsto
p_+(\hatgamma):={\textstyle\bigsqcup\limits_{(x,\myp\bar{y})\in\myp
\hat{\gamma\myp}}} \{(x,\bar{y}+x)\}\in \varGamma_{\mathcal{Z}}^{\sharp},\\
\label{eq:pY} \varGamma_{\mathcal{Z}}^{\sharp}\ni \hatgamma&\mapsto
p_{\mathfrak{X}}(\hatgamma):={\textstyle\bigsqcup\limits_{(x,\myp\bar{y})\in\myp
\hat{\gamma\myp}}}\{\bar{y}\}\in \varGamma_{\mathfrak{X}}^{\sharp},\\
\label{eq:pr2} \varGamma_{\frakX}^\sharp\ni\bar{\gamma}&\mapsto
\mathfrak{p}(\bar{\gamma}):={\textstyle\bigsqcup\limits_{\bar{y}\in
\bar{\gamma}}}\mathfrak{p}(\bar{y})\in\varGamma_{X}^\sharp.
\end{align}
To verify that the map
$p_+:\varGamma_{\mathcal{Z}}^\sharp\to\varGamma_{\mathcal{Z}}^\sharp$
is measurable, note that for a cylinder set
$$
C_{B_1\times \bar{B}}^{\myp n}
=\{\hatgamma\in\varGamma_{\mathcal{Z}}: \,\hatgamma(B_1\times
\bar{B})=n\}\in \mathcal{B}(\varGamma_{\mathcal{Z}}^\sharp),
$$
with $B_1\in\mathcal{B}(X)$, $\bar{B}\in\mathcal{B}(\mathfrak{X})$
and $n\in\ZZ_+$, its pre-image under $p_+$ is given by
$$
p_+^{-1}(C_{B_1\times \bar{B}}^{\myp n})=C_{A}^{\myp
n}=\{\hatgamma\in\varGamma_{\mathcal{Z}}^\sharp: \,\hatgamma(A)=n\},
$$
where
$$
A:=\{(x,\bar{y})\in\mathcal{Z}:(x,\bar{y}+x)\in
B_1\times\bar{B}\}\in\mathcal{B}(\mathcal{Z}),
$$
since the indicator
$\mathbf{1}_A(x,\bar{y})=\mathbf{1}_{B_1}(x)\cdot
\mathbf{1}_{\bar{B}}(\bar{y}+x)$ is obviously a measurable function
on $\mathcal{Z}=X\times \mathfrak{X}$. Similarly, for
$p_{\mathfrak{X}}:\varGamma_{\mathcal{Z}}^\sharp\to\varGamma_{\mathfrak{X}}^\sharp$
(see (\ref{eq:pY})) we have
$$
p_{\mathfrak{X}}^{-1}(C_{\bar{B}}^{\myp n})=C_{X\times
\bar{B}}^{\myp n}=\{\hatgamma\in\varGamma_{\mathcal{Z}}^\sharp:
\,\hatgamma(X\times\bar{B})=n\}\in
\mathcal{B}(\varGamma_{\mathcal{Z}}^\sharp),
$$
since $X\times\bar{B}\in \mathcal{B}(\mathcal{Z})$. Finally,
measurability of the projection
$\mathfrak{p}:\varGamma_{\mathfrak{X}}^\sharp\to\varGamma_X^\sharp$
(see (\ref{eq:pr2})) was shown in \cite[Section 3.3, p.~455]{BD3}.
As a result, the composition (\ref{eq:circ}) is measurable, as
claimed.
\endproof

Let us define a measure on $\varGamma^\sharp_{X}$ as the
push-forward of $\hat{\mu}$ (see Definition~\ref{def:muhat}) under
the map $\mathfrak{q}$ defined in (\ref{proj1}), (\ref{eq:proj}):
\begin{equation}\label{eq:mu*}
\mathfrak{q}^{\ast}\hat{\mu}(A)
\equiv\hat{\mu}\myp(\mathfrak{q}^{-1}(A)),\qquad
A\in\mathcal{B}(\varGamma_{X}^\sharp),
\end{equation}
or equivalently
\begin{equation}\label{eq:mu*1}
\int_{\varGamma_{X}^\sharp}F(\gamma
)\,\mathfrak{q}^{\ast}\hat{\mu}(\rd\gamma)
=\int_{\varGamma^\sharp_{\mathcal{Z}}} \mynn
F(\mathfrak{q}(\hatgamma)) \,\hat{\mu}({\rd}\hatgamma),\qquad
F\in\mathrm{M}_+(\varGamma_{X}^\sharp).
\end{equation}
The next general result shows that this measure may be identified
with the original cluster measure $\mu_{\rm cl}$.
\begin{theorem}\label{th:mucl}
The measure \textup{(\ref{eq:mu*})} coincides with the cluster
measure $\mu_{\rm cl}$,
\begin{equation}\label{eq:cl*}
\mu_{\rm cl}=\mathfrak{q}^{\ast}\hat{\mu}\equiv \hat{\mu}\circ
\mathfrak{q}^{-1}.
\end{equation}
\end{theorem}

\proof Let us evaluate the Laplace transform of the measure
$\mathfrak{q}^{\ast}\hat{\mu}$. For any function
$f\in\mathrm{M}_+(X)$, we obtain, using (\ref{eq:mu*1}),
(\ref{eq:proj}) and (\ref{eq:muhat}),
\begin{align*}
L_{\mathfrak{q}^{\ast}\myn\hat{\mu}}(f)&= \int_{\varGamma
_{X}^\sharp}\exp(-\langle
f,\gamma\rangle)\,\mathfrak{q}^{\ast}\hat{\mu}(\rd\gamma)
= \int_{\varGamma_{\mathcal{Z}}^\sharp}\exp(-\langle
f,\mathfrak{q}(\hatgamma)\rangle)\,\hat{\mu}(\rd\hatgamma)\\[.1pc]
&=\int_{\varGamma_{X}}\biggl(\int_{\varGamma_{\mathfrak{X}}^\sharp}
\exp\biggl(-\sum_{x\in\gamma_{\rm c}}
f(\mathfrak{p}(\bar{y}_x+x))\biggr)
{\textstyle\bigotimes\limits_{x\in\gamma_{\rm c}}}\eta(\rd\bar{y}_x)\biggr)
\,\mu_{\rm c}(\rd\gamma_{\rm c})\\[.1pc]
&=\int_{\varGamma_
X^\sharp}\biggl(\int_{\varGamma_{\mathfrak{X}}^\sharp}\,
\prod_{x\in\gamma_{\rm c}}\exp\bigl(-f(\bar{y}_x+x)\bigr)
{\textstyle\bigotimes\limits_{x\in\gamma_{\rm c}}}\eta(\rd\bar{y}_x)\biggr)\,
\mu_{\rm c}(\rd\gamma_{\rm c})\\[.1pc]
&=\int_{\varGamma_X^\sharp} \prod_{x\in\gamma_{\rm
c}}\biggl(\int_{\mathfrak{X}} \exp \biggl(- \sum_{y\in
\bar{y}}f(y+x)\biggr) \,\eta(\rd\bar{y})\biggr)\,\mu_{\rm
c}(\rd\gamma_{\rm c}),
\end{align*}
which coincides with the Laplace transform (\ref{laplace-G}) of the
cluster measure $\mu_{\rm cl}$.
\endproof

\subsection{Gibbs cluster measure via an auxiliary Gibbs measure}\label{sec:2.4}
In this paper, we are concerned with \emph{Gibbs cluster point
processes}, for which the
distribution of cluster centers is
given by some
Gibbs measure $\g\in\mathscr{G}(\theta,\U\myp)$ on the proper
configuration space $\varGamma_{X}$ (see the Appendix), specified by
a
\emph{reference measure} $\theta$ on $X$ and an \textit{interaction
potential} \,$\U:\varGamma^{\myp 0}_{X}\rightarrow \RR\cup
\{+\infty\}$, where $\varGamma_{X}^{\myp0}\subset\varGamma_X$ is the
subspace of finite configurations in $X$. We assume that the set
\mypp$\mathscr{G}(\theta,\varPhi)$ of all Gibbs measures on
$\varGamma_{X}$ associated with $\theta$ and $\U$ is
non-empty.\footnote{For various sufficient conditions, consult
\cite{Preston,Ruelle}; see also references in the Appendix.}
%

Specializing Definition \ref{def:muhat} to the Gibbs case, the
corresponding auxiliary measure $\hatg$ on the (proper)
configuration space $\varGamma_\mathcal{Z}$ is given by (cf.\
(\ref{eq:muhat0}), (\ref{eq:muhat}))
\begin{equation}\label{eq:ghat0}
\hatg(\rd\hatgamma)=\g(p_{X}(\rd\hatgamma))
\,{\textstyle\bigotimes\limits_{z\in\hat{\gamma\myp}}}
\,\eta(p_\frakX(\rd{z})),\qquad \hatgamma\in\varGamma_{\mathcal{Z}},
\end{equation}
or equivalently, for any function
$F\in{\mathrm{M}}_+(\varGamma_{\mathcal{Z}})$,
\begin{equation}\label{eq:ghat}
\int_{\varGamma_{\mathcal{Z}}} F(\hatgamma)\,\hatg(\rd\hatgamma)
=\int_{\varGamma_{X}}\!\biggl(\int_{{\mathfrak{X}}^\infty}
F\Bigl(\myp{\textstyle\bigcup\nolimits_{x\in\gamma}}\{(x,\bar{y})\}
\myn\Bigr)\,{\textstyle\bigotimes\limits_{x\in\gamma}}
\,\eta(\rd\bar{y})\mynn\biggr)\,\g(\rd\gamma).
\end{equation}

\begin{remark} Vector $\bar{y}$ in each pair $z=(x,\bar{y})\in
\mathcal{Z}$ may be interpreted as a \emph{mark} attached to the
point $x\in X$, so that $\hatgamma$ becomes a marked configuration,
with the mark space $\mathfrak{X}$ (see \cite{DVJ1,KunaPhD,KKS98}).
The corresponding \textit{marked configuration space} is defined by
\begin{equation}\label{eq:Gamma(X)}
\varGamma_{X}(\mathfrak{X}):=\{\hatgamma\in\varGamma_{\mathcal{Z}}:
p_{X}(\hatgamma)\in\varGamma_X\}\subset \varGamma_{\mathcal{Z}}.
\end{equation}
In other words, $\varGamma_{X}(\mathfrak{X})$ is the set of
configurations in $\varGamma_{\mathcal{Z}}$ such that the collection
of their $x$-coordinates is a (proper) configuration in
$\varGamma_{X}$. Clearly, $\varGamma_{X}(\mathfrak{X})$ is a Borel
subset of $\varGamma_{\mathcal{Z}}$, that is,
$\varGamma_{X}(\mathfrak{X})\in\mathcal{B}(\varGamma_{\mathcal{Z}})$.
Since $\g(\varGamma_X)=1$, we have
$\hatg(\varGamma_{X}(\mathfrak{X}))=1$.
\end{remark}

Finally, owing to the general Theorem \ref{th:mucl} (see
(\ref{eq:cl*})), the corresponding Gibbs cluster measure $\gcl$ on
the configuration space $\varGamma_X$ is represented as a
push-forward of the measure $\hatg$ on $\varGamma_{\mathcal{Z}}$
under the map $\mathfrak{q}$ defined in (\ref{proj1}),
(\ref{eq:proj}):
\begin{equation}\label{eq:gcl*}
\gcl=\mathfrak{q}^*\hatg\equiv \hatg\circ \mathfrak{q}^{-1}.
\end{equation}

Our next goal is to show that $\hatg$ is a \textit{Gibbs measure} on
$\varGamma_{\mathcal{Z}}$, with the reference measure $\sigma$
defined as a product measure on the space
$\mathcal{Z}=X\times\mathfrak{X}$,
\begin{equation}\label{eq:sigma}
\sigma:=\theta \otimes\eta,
\end{equation}
and with the interaction potential \mypp$\hat{\U}:\varGamma^{\myp
0}_{\mathcal{Z}}\to\RR\cup\{+\infty\}$ given by
\begin{equation}\label{eq:hatU}
\hat{\U}(\hatgamma):=\left\{\begin{array}{ll}
\U(p_{X}(\hatgamma)),\quad&
\hatgamma\in\varGamma^{\myp 0}_{\mathcal{Z}}\cap \varGamma_{X}(\mathfrak{X}),\\[.3pc]
+\infty,&\hatgamma\in\varGamma^{\myp
0}_{\mathcal{Z}}\setminus\varGamma_{X}(\mathfrak{X}),
\end{array}\right.
\end{equation}
where $p_{X}$ is the projection defined in (\ref{eq:pX}). The
corresponding functionals of energy $\hat{E}(\hat{\xi}\myp)$ and
interaction energy $\hat{E}(\hat{\xi},\hatgamma)$ ($\hat{\xi} \in
\varGamma_{\mathcal{Z}}^{\myp0}$, $\hatgamma\in
\varGamma_{\mathcal{Z}}$) are then given by (see (\ref{eq:E}) and
(\ref{eq:E-E}))
\begin{gather}
\label{eq:Ehat}
\hat{E}(\hat{\xi}\myp):=\sum_{\hat{\xi}'\subset
\myp\hat{\xi}}\hat{\varPhi}(\hat{\xi}'),\\
\label{eq:E-Ehat} \hat{E}(\hat{\xi},\hatgamma):=\left\{
\begin{array}{ll}
\displaystyle\sum_{\hatgamma\supset\myp\hatgamma'\myn\in\varGamma_{\mathcal{Z}}^{\myp0}}
\!\hat{\varPhi}(\hat{\xi}\cup\hatgamma'), & \displaystyle
\sum_{\hatgamma\supset\myp\hatgamma'\myn\in\varGamma_{\mathcal{Z}}^{\myp0}}
\!|\myp\hat{\varPhi}(\hat{\xi}\cup\hatgamma')| <\infty , \\[1.7pc]
\displaystyle\ \ +\infty  & \quad \text{otherwise}.
\end{array}
\right.
\end{gather}

The following ``projection''property of the energy is obvious from
the definition (\ref{eq:hatU}) of the potential
$\hat{\U}$.
\begin{lemma}\label{lm:invariance}
For any configurations\/ $\hat{\xi}\in
\varGamma_{\mathcal{Z}}^{\myp0}$ and $\hatgamma\in
\varGamma_{\mathcal{Z}}$, we have
\begin{equation*}
\hat{E}(\hat{\xi}\myp)=E(p_{X}(\hat{\xi}\myp)),\qquad
\hat{E}(\hat{\xi},\hatgamma)=E(p_{X}(\hat{\xi}\myp),p_{X}(\hatgamma)).
\end{equation*}
\end{lemma}

%

\begin{theorem}\label{th:g-hat}
\textup{(a)} \,Let\/ $\g\in{\mathscr{G}}(\theta,\U\myp)$ be a Gibbs
measure on the configuration space\/ $\varGamma_X$, and let\/
$\hatg$ be the corresponding probability measure on the
configuration space\/ $\varGamma_{\mathcal{Z}}$ \textup{(}see
\textup{(\ref{eq:ghat0})} or \textup{(\ref{eq:ghat})}\textup{)}.
Then\/ $\hatg\in {\mathscr{G}}(\sigma,\hat{\U}\myp)$, i.e., $\hatg$
is a Gibbs measure on\/ $\varGamma_{\mathcal{Z}}$ with the reference
measure\/ $\sigma$ and the interaction potential\/ $\hat{\U}$
defined by \textup{(\ref{eq:sigma})} and \textup{(\ref{eq:hatU})},
respectively.

\textup{(b)} \,If the measure\/ $\g\in{\mathscr{G}}(\theta,\U\myp)$
has a finite correlation function\/ $\kappa_{\myn\g}^{n}$ of some
order\/ $n\in\NN$ \textup{(}see the definition
\textup{(\ref{corr-funct})} in the Appendix\textup{)}, then the
correlation function\/ $\kappa_{\myn\hatg}^{n}$ of the measure\/
$\hatg\in {\mathscr{G}}(\sigma,\hat{\U}\myp)$ is given by
\begin{equation}\label{corr-funct0}
\kappa_{\myn\hatg}^{n}(z_{1},\dots,z_{n})
=\kappa_{\myn\g}^{n}\bigl(p_{X}(z_{1}),\dots,p_{X}(z_{n})\bigr),\qquad
z_1,\dots,z_n\in \mathcal{Z}.
\end{equation}
\end{theorem}

\proof (a)
\,In order to show that $\hatg\in
{\mathscr{G}}(\sigma,\hat{\U}\myp)$, it suffices to check that
$\hatg$ satisfies Nguyen--Zessin's equation on
$\varGamma_{\mathcal{Z}}$ (see equation (\ref{eq:NZ}) in the
Appendix), that is, for any non-negative, $\mathcal{B}(\mathcal{Z})
\times \mathcal{B}(\varGamma_{\mathcal{Z}})$-measurable function
$H(z,\hatgamma)$ it holds
\begin{equation}\label{eq:NZ-g-hat-exp}
\int_{\varGamma_{\mathcal{Z}}} \sum_{z\in\hat{\gamma\myp}}
\!H(z,\hatgamma)\,\hatg(\rd\hatgamma)
=\int_{\varGamma_{\mathcal{Z}}}\!
\left(\int_{\mathcal{Z}}\!H(z,\hatgamma\cup
\{z\})\,\re^{-\hat{E}(\{z\},\myp\hatgamma )}\,\sigma
(\rd{z})\mynn\right)\hatg(\rd\hatgamma).
\end{equation}
Using the disintegration formula (\ref{eq:muhat}), the left-hand
side of (\ref{eq:NZ-g-hat-exp}) can be represented as
\begin{multline}\label{eq:cond-eta}
\int_{\varGamma_X}\!\biggl(\int_{{\frakX}^\infty} \,\sum_{x\in
\gamma}
H\bigl(x,\bar{y}_x;\,{\textstyle\bigcup\limits_{x'\in\gamma}}\{({\textstyle
x'}\mynn,\bar{y}_{x'\mynn})\}\bigr)
{\textstyle\bigotimes\limits_{x'\in\gamma}}\eta(\rd\bar{y}_{x'\myn})\mynn\biggr)\,\g(\rd\gamma) \\
=\int_{\varGamma_X}\!\left(\mynn\sum_{x\in\gamma}
\int_{\varGamma_{\frakX}} \mathbf{1}_{\gamma}(x)\,
H\bigl(x,\bar{y}_x;\,{\textstyle\bigcup\limits_{x'\in\gamma}}\{({\textstyle
x'}\mynn,\bar{y}_{x'\mynn})\}\bigr)
{\textstyle\bigotimes\limits_{x'\in\gamma}}\,\eta(\rd\bar{y}_{x'\mynn})\mynn\right)\g(\rd\gamma).
\end{multline}
Applying Nguyen--Zessin's equation to the Gibbs measure $\g$ with
the function
\begin{equation*}
H_0(x,\gamma):=\int_{\varGamma_{\frakX}} \mathbf{1}_{\gamma}(x)\,
H\Bigl(x,\bar{y}_x;\,{\textstyle\bigcup\limits_{x'\in\gamma}}\{({\textstyle
x'}\mynn,\bar{y}_{x'\mynn})\}\Bigr)
{\textstyle\bigotimes\limits_{x'\in\gamma}}\,\eta(\rd\bar{y}_{x'\mynn}),
\end{equation*}
we see that the right-hand side of (\ref{eq:cond-eta}) takes the
form
\begin{equation}\label{eq:E_g}
\int_{\varGamma_X}\!\biggl(\mynn\sum_{x\in \gamma}
H_0(x,\gamma)\mynn\biggr)\,\g(\rd\gamma)= \int_{\varGamma_X}
\!\biggl(\int_{X} H_0(x,\gamma\cup
\{x\})\,\re^{-E(\{x\},\myp\gamma)}\,\theta(\rd
x)\mynn\biggr)\,\g(\rd\gamma).
\end{equation}
Similarly,
exploiting the product structure of the measures
\,$\sigma=\theta\otimes\eta$ and
$$
{\textstyle\bigotimes\limits_{x'\in\gamma\cup\{x\}}}
\mynn\eta(\rd\bar{y}_{x'\mynn})=
{\textstyle\bigotimes\limits_{x'\in\gamma}}\,\eta(\rd\bar{y}_{x'\mynn})
\otimes \eta(\rd\bar{y}_{x}),
$$
and using Lemma \ref{lm:invariance}, the right-hand side of
(\ref{eq:NZ-g-hat-exp}) is reduced to
\begin{align*}
&\int_{\varGamma_X}\!\biggl(\int_{\mathfrak{X}^{\infty}}
\biggl(\int_{X\times\mathfrak{X}}
H\bigl(x,\bar{y}_x;{\textstyle\bigcup\limits_{x'\in\gamma}}\{({\textstyle
x'}\mynn,\bar{y}_{x'\mynn})\}\cup\{(x,\bar{y}_x)\}\bigr)\\
&\hspace{9.97pc}
\times\re^{-E(\{x\},\myp\gamma)}\,\eta(\rd\bar{y}_x)\,\theta(\rd
x)\biggr) {\textstyle\bigotimes\limits_{x'\in\gamma}}
\,\eta(\rd\bar{y}_{x'\mynn})\biggr)\,\g(\rd\gamma)\\
&= \int_{\varGamma_X}\!\biggl(\int_X
\biggl(\int_{\mathfrak{X}^{\infty}}
\mathbf{1}_{\gamma}(x)\,H\bigl(x,\bar{y}_x;{\textstyle\bigcup\limits_{x'\in\gamma}}
\{({x'}\mynn,\bar{y}_{x'\mynn})\}\bigr)\\
&\hspace{9.97pc}\times
\re^{-E(\{x\},\myp\gamma)}{\textstyle\bigotimes\limits_{x'\in\gamma}}\,\eta(\rd\bar{y}_{x'\mynn})
\biggr)\,\theta(\rd x)\biggr)\,\g(\rd\gamma)\\
&=\int_{\varGamma_X}\!\biggl(\int_{X} H_0(x,\gamma\cup
\{x\})\,\re^{-E(\{x\},\myp\gamma)}\,\theta(\rd{x})\biggr)\,\g(\rd\gamma),
\end{align*}
thus coinciding with (\ref{eq:E_g}). This proves equation
(\ref{eq:NZ-g-hat-exp}), hence $\hatg\in
{\mathscr{G}}(\sigma,\hat{\U}\myp)$.

(b) \,Let $f\in C_0(\mathcal{Z}^{n})$ be a symmetric function.
According to the disintegration formula (\ref{eq:muhat}) applied to
the function
\begin{equation*}
F(\hatgamma):=\sum_{\{z_{1}\myn,\dots,\myp z_{n}\} \subset
\hatgamma}f(z_{1},\dots,z_{n}),
\end{equation*}
we have
\begin{equation}\label{eq:sum-phi}
\int_{\varGamma_{\mathcal{Z}}} \mynn
F(\hatgamma)\,\hatg(\rd\hatgamma)=
\int_{\varGamma_{X}}\!\myn{}\sum_{\{x_{1}\myn,\dots,\myp
x_{n}\}\subset\gamma}\mynn \phi (x_{1},\dots,x_{n})\, \g(\rd\gamma),
\end{equation}
where
\begin{equation*}
\phi (x_{1},\dots,x_{n}):=\int_{\mathfrak{X}^{n}} \mynn
f((x_{1},\bar{y}_{1}\myn),\dots,(x_{n},\bar{y}_{n}))\,{\textstyle\bigotimes\limits_{i=1}^n}\,
\eta(\rd\bar{y}_{i})\in C_0(X^n).
\end{equation*}
Applying the definition of the correlation function
$\kappa_{\myn\g}^n$ (see (\ref{corr-funct})) and using that
$\theta(\rd{x})\otimes\eta(\rd\bar{y})=\sigma(\rd{x}\times\rd\bar{y})$,
we obtain from (\ref{eq:sum-phi})
\begin{equation*}
\int_{\varGamma_{\mathcal{Z}}} \mynn
F(\hatgamma)\,\hatg(\rd\hatgamma)
=\frac{1}{n!}\int_{\mathcal{Z}^{n}}f(z_{1},\dots,z_{n})
\,\kappa_{\myn\g}^{n}\bigl(p_{X}(z_{1}),\dots,p_{X}(z_{n})\bigr)
\,{\textstyle\bigotimes\limits_{i=1}^n}\,\sigma(\rd{z}_{i}),
\end{equation*}
and equality (\ref{corr-funct0}) follows.
\endproof

In the rest of this subsection, $\GR$ denotes the subclass of Gibbs
measures in $\mathscr{G}$ (with a given reference measure and
interaction potential) that satisfy the so-called \textit{Ruelle
bound} (see the Appendix, formula (\ref{eq:RB})).

\begin{corollary}\label{cor:GR}
We have\/ $\g\in\GR(\theta,\U\myp)$ if and only if\/ $\hatg\in
\GR(\sigma,\hat{\U}\myp)$,
\end{corollary}

\proof Follows directly from formula (\ref{corr-funct0}).
\endproof

The following statement is, in a sense, converse to Theorem
\ref{th:g-hat}\myp(a).
\begin{theorem}\label{th:converse}
If\/ $\varpi \in {\mathscr{G}}(\sigma ,\hat{\U}\myp)$ then
$\g:=p_{X}^{\ast}\varpi \in {\mathscr{G}}(\theta,\U\myp)$. Moreover,
if\/ $\g\in\GR(\theta,\U\myp)$ then\/ $\varpi =\hatg$.
\end{theorem}

\proof Applying Nguyen--Zessin's equation (\ref{eq:NZ}) to the
measure $\varpi$ and using the cylinder structure of the interaction
potential, we have
\begin{align*}
\int_{\varGamma_{X}}&\sum_{x\in \gamma}H(x,\gamma )\,p_{X}^{\ast}
\varpi(\rd\gamma)=\int_{\varGamma_{\mathcal{Z}}}\sum_{x\in
p_{X}\hatgamma} H(x,p_{X}\hatgamma)\,\varpi(\rd\hatgamma) \\
&=\int_{\varGamma_{\mathcal{Z}}}\!\left(
\int_{\mathcal{Z}}H(p_{X}z,p_{X}(\hatgamma\cup \{z\}
)\,\re^{-E(\{p_{X}z\},\,p_{X}\hatgamma)}\,\theta\otimes
\eta(\rd{z})\mynn\right) \varpi (\rd\hatgamma) \\
&=\int_{\varGamma_{\mathcal{Z}}}\!\left(
\int_{X}H(x,p_{X}\hatgamma\cup
\{x\})\,\re^{-E(\{x\},\,p_{X}\hatgamma)}\,\theta(\rd x)\mynn\right)
\varpi(\rd\hatgamma) \\
&=\int_{\varGamma_{X}}\!\left(\int_{X} H(x,\gamma
\cup\{x\})\,\re^{-E(\{x\},\myp\gamma)}\,\theta(\rd x)\right)
p_{X}^{\ast} \varpi(\rd\gamma).
\end{align*}
Thus, the measure $p_{X}^*\varpi$ satisfies Nguyen--Zessin's
equation and so, by Theorem \ref{th:NZR}, belongs to the Gibbs class
$\mathscr{G}(\theta,\varPhi)$.

Next, in order to prove that $\varpi =\hatg$, by Proposition
\ref{pr:k=k} it suffices to show that the measures $\varpi$ and
$\hatg$ have the same correlation functions. Note that the
correlation function $\kappa_{\varpi}^{n}$ can be written in the
form \cite[\S\myp2.3, Lemma~2.3.8]{KunaPhD}
\begin{align*}
\kappa_{\varpi}^{n}(z_1,\dots,z_n) &=\re^{-\hat{E}(\{z_1,\dots,\myp
z_n\})}\int_{\varGamma_{\mathcal{Z}}}
\re^{-\hat{E}(\{z_1,\dots,\myp z_n\},\myp\hatgamma)}\,\varpi(\rd\hatgamma)\\
&=\re^{-E(\{p_{X}(z_1),\dots,\mypp
p_{X}(z_n)\})}\int_{\varGamma_{X}}
\re^{-E(\{p_{X}(z_1),\dots,\mypp p_{X}(z_n)\},\myp\gamma)}\,p_{X}^{\ast }\varpi(\rd\gamma)\\
&= \re^{-E(\{p_{X}(z_1),\dots,\mypp p_{X}
(z_n)\})}\int_{\varGamma_{X}}
\re^{-E(\{p_{X}(z_1),\dots,\mypp p_{X}(z_n)\},\myp\gamma)}\,\g(\rd\gamma) \\
&=\kappa_{\myn\g}^{n}\bigl(p_{X}(z_1),\dots, p_{X}(z_n)\bigr).
\end{align*}
Therefore, on account of Theorem \ref{th:g-hat}\myp(b) we get $
\kappa_{\varpi}^{n}(z_1,\dots,z_n)=\kappa_{\myn\hatg}^{n}(z_1,\dots,z_n)$
for all $z_1,\dots,z_n\in\mathcal{Z}$ \,($z_i\ne z_j$), as required.
\endproof

In the next corollary, $\ext\mathscr{G}$ denotes the set of
\textit{extreme points} of the class $\mathscr{G}$ of Gibbs measures
with the corresponding reference measure and interaction potential
(see the Appendix).
\begin{corollary}\label{cor:2.9}
Suppose that\/ $\g\in\GR(\theta,\U\myp)$. Then\/
$\g\in\ext\mathscr{G}(\theta,\U\myp)$  if and only if\/
$\hatg\in\ext\mathscr{G}(\sigma,\hat{\U}\myp)$.
\end{corollary}
\proof Let $\g\in \GR(\theta,\U\myp)\cap
\ext\mathscr{G}(\theta,\U\myp)$. Assume that
$\hatg=\frac12\myp(\mu_{1}+\mu_{2})$ with some $\mu_{1},\mu_{2}\in
\mathscr{G}(\sigma,\hat{\U}\myp)$. Then
$\g=\frac12\myp(\g_{1}+\g_{2})$, where $\g_{i}=p_{X}^{\ast}
\mu_{i}\in\mathscr{G}(\theta,\U\myp)$. Since
$\g\in\ext\mathscr{G}(\theta,\U\myp)$, this implies that
$\g_{1}=\g_{2}=\g$. In particular,
$\g_{1},\g_{2}\in\GR(\theta,\U\myp)$ and by Theorem
\ref{th:converse} we obtain that
$\mu_{1}=\hatg_{1}=\hatg=\hatg_{2}=\mu _{2}$, which implies
$\hatg\in \ext \mathscr{G}(\sigma,\hat{\U}\myp)$.

Conversely, let $\hatg\in\ext\mathscr{G}(\sigma,\hat{\U}\myp)$ and
$\g=\frac12\myp(\g_{1}+\g_{2})$ with $\g_{1},\g_{2}\in
\mathscr{G}(\theta,\U\myp)$. Then
$\hatg=\frac12\myp(\hatg_{1}+\hatg_{2})$, hence
$\hatg_{1}=\hatg_{2}=\hatg\in\GR(\sigma,\hat{\U}\myp)$, which
implies by Theorem \ref{th:converse} that
$\g_{1}=p_{X}^{\ast}\hatg_{1}=p_{X}^{\ast}\hatg_{2}=\g_{2}$. Thus,
$\g\in\ext\mathscr{G}(\theta,\U\myp)$.
\endproof

\subsection{Criteria of local finiteness and simplicity of
the Gibbs cluster process}\label{sec:2.5}

Let us give conditions sufficient for the Gibbs CPP to be (a)
locally finite, and (b) simple. For a given Borel set
$B\in\calB(X)$, consider a set-valued function (referred to as the
\textit{droplet cluster})
\begin{equation}\label{eq:D}
D_B(\bar{y}):={\textstyle\bigcup\limits_{y_i\in\bar{y}}}
(B-y_i),\qquad \bar{y}\in\frakX.
\end{equation}
Let us also denote by $N_B(\bar{y})$ the number of coordinates of
the vector $\bar{y}=(y_i)$ falling in the set $B\in\mathcal{B}(X)$,
\begin{equation}\label{eq:N(y)}
N_B(\bar{y}):=\sum_{y_i\in\bar{y}} \mathbf{1}_B(y_i),\qquad
\bar{y}\in\frakX,
\end{equation}
In particular, $N_X(\bar{y})$ is the ``dimension'' of $\bar{y}$,
that is, the total number of its coordinates (recall that
$\bar{y}\in\frakX=\bigsqcup_{\myp n=0}^{\myp\infty} X^n$, see
(\ref{eq:calX})).

\begin{theorem}\label{th:properClusterGibbs} Let
$\gcl$ be a Gibbs cluster measure on the generalized configuration
space $\varGamma_ X^\sharp$.

\textup{(a)} \,Assume that the
correlation function $\kappa_{\myn\g}^1$ of the measure
$\g\in\mathscr{G}(\theta,\U\myp)$ is bounded.
Then, in order that $\gcl$-a.a.\
configurations $\gamma\in \varGamma_{X}^\sharp$ be locally finite,
it is sufficient that the following two conditions hold\textup{:}

{\rm (a-i)} \,for any compact set\/ $B\in{\mathcal{B}}(X)$, the
number of coordinates of the vector\/ $\bar{y}\in\frakX$ in $B$ is
a.s.-finite,
\begin{equation}\label{eq:condA1}
N_B(\bar{y})<\infty\quad \text{for} \,\,\eta\text{-a.a.} \
\bar{y}\in\frakX;
\end{equation}

{\rm (a-ii)} \,for any compact set $B\in{\mathcal{B}}(X)$, the mean
$\theta$-measure of the droplet cluster\/ $D_B(\bar{y})$ is finite,
\begin{equation}\label{eq:condA2}
\int_{\frakX} \theta(D_B(\bar{y}))\, \eta(\rd\bar{y})<\infty\myp.
\end{equation}

\textup{(b)} \,In order that $\gcl$-a.a.\ configurations $\gamma\in
\varGamma_{X}^\sharp$ be simple, it is sufficient that the following
two conditions hold\textup{:}

{\rm (b-i)} \,for any $x\in X$, vector $\bar{y}$ contains a.s.\ no
more than one coordinate $y_i=x$,
\begin{equation}\label{eq:condB1}
\sup_{x\in X} N_{\{x\}}(\bar{y})\le 1\quad \text{for}
\,\,\eta\text{-a.a.} \ \bar{y}\in\frakX;
\end{equation}

{\rm (b-ii)} \,for any $x\in X$, the ``point'' droplet cluster
$D_{\{x\}}(\bar{y})$
has a.s.\ zero $\theta$-measure,
\begin{equation}\label{eq:condB2}
\theta\bigl(D_{\{x\}}(\bar{y})\bigr)=0\quad \text{for}
\,\,\eta\text{-a.a.} \ \bar{y}\in\frakX.
\end{equation}
\end{theorem}

For the proof of part (a) of this theorem, we need a reformulation
(stated as Proposition \ref{prop1} below) of the condition (a-ii),
which will also play an important role in utilizing the projection
construction of the Gibbs cluster measure (see Section~\ref{sec:3}
below). For any Borel subset $B\in\mathcal{B}(X)$, denote
\begin{equation}\label{eq:K}
\mathcal{Z}_{B}:=\{z\in \mathcal{Z}:\,\mathfrak{q}(z)\cap
B\ne\emptyset \}\in\mathcal{B}(\mathcal{Z}),
\end{equation}
where $\mathfrak{q}(z)=\bigsqcup_{y_i\in
p_\frakX(z)}\{y_i+p_{X}(z)\}$ (see (\ref{proj1})). That is to say,
the set $\mathcal{Z}_{B}$ consists of all points $z=(x,\bar{y})\in
\mathcal{Z}$ such that, under the ``projection'' $\mathfrak{q}$ onto
the space $X$, at least one coordinate $y_i+x$ ($y_i\in\bar{y}$)
belongs to the set $B\subset X$.
\begin{proposition}\label{prop1}
For any $B\in\mathcal{B}(X)$,
the condition \mbox{\textup{(a-ii)}} of\/ Theorem
\textup{\ref{th:properClusterGibbs}\myp(a)} is necessary and
sufficient in order that\/ $\sigma(\mathcal{Z}_{B})<\infty$, where\/
$\sigma=\theta\otimes\eta$.
\end{proposition}

\proof[Proof\/ of\/ Proposition \textup{\ref{prop1}}]
By definition (\ref{eq:K}), $(x,\bar{y})\in \mathcal{Z}_B$ if and
only if $x\in \bigcup_{y_i\in\bar y} (B-y_i)\equiv D_B(\bar y)$ (see
(\ref{eq:D})). Hence,
\begin{align*}
\sigma(\mathcal{Z}_{B})&=\int_{\frakX}\!\left(\int_ X
\mathbf{1}_{D_{B}(\bar y)}(x)\,\theta(\rd x)\right)\eta(\rd\bar{y})
=\int_{\frakX}\theta\bigl(D_B(\bar y)\bigr)\,\eta(\rd\bar{y}),
\end{align*}
and we see that the bound $\sigma(\mathcal{Z}_{B})<\infty$ is
nothing else but condition (\ref{eq:condA2}).
\endproof

\proof[Proof\/ of\/ Theorem \textup{\ref{th:properClusterGibbs}}]
(a) \,Let $B\subset X$ be a compact set. By Proposition \ref{prop1},
condition (a-ii) is equivalent to $\sigma(\mathcal{Z}_B)<\infty$. On
the other hand, by Theorem \ref{th:g-hat}\myp(b) we have
$\kappa_{\myn\hatg}^1(x,\bar{y})=\kappa_{\myn\g}^1(x)$. Hence,
$\kappa_{\myn\hatg}^1$ is bounded, and by Remark
\ref{rm:kappa<const} (see the Appendix) it follows that
$\hatgamma(\mathcal{Z}_B)<\infty$ ($\hatg$-a.s.). According to the
projection representation $\gcl=\mathfrak{q}^*\hatg$ (see
(\ref{eq:gcl*})) and in view of condition \mbox{(a-i)}, this implies
that, almost surely, a projected configuration
$\gamma=\mathfrak{q}(\hatgamma)=\bigsqcup_{z\in
\hat{\gamma\myp}}\mathfrak{q}(z)$ contributes no more than
finitely many points to the set $B\subset
\mathfrak{q}(\mathcal{Z}_B)$, that is, $\gamma(B)<\infty$
($\gcl$-a.s.), which completes the proof of part~(a).


(b) \,It suffices to prove that, for any compact set
$\varLambda\subset X$,
there are $\gcl$-a.s.\ no cross-ties between the clusters whose
centers belong to $\varLambda$.
That is, we must show that $\gcl(A_\varLambda)=0$, where the set
$A_\varLambda\in\mathcal{B}(\varGamma_X\times \frakX^2)$ is defined
by
\begin{equation}\label{eq:A-Lambda}
A_\varLambda:=\{(\gamma,\bar{y}_1,\bar{y}_2): \,\exists\,
x_1,x_2\in{}\gamma\cap\varLambda,\ \exists\,y_1\in\bar{y}_1:
\,x_1+y_1-x_2\in\bar{y}_2\},
\end{equation}
Applying the disintegration formula (\ref{eq:ghat}), we obtain
\begin{align}\label{eq:Psi3}
\gcl(A_\varLambda)&=\int_{\varGamma_X} \mynn
F(\gamma)\,\g(\rd\gamma),
\end{align}
where
\begin{equation}\label{eq:F}
F(\gamma):=\int_{\frakX^2}
\mathbf{1}_{A_\varLambda}(\gamma,\bar{y}_1,\bar{y}_2)\,
\eta(\rd\bar{y}_1)\,\eta(\rd\bar{y}_2),\qquad \gamma\in\varGamma_X.
\end{equation}
Note that, according to the definition (\ref{eq:A-Lambda}),
$F(\gamma)\equiv F(\gamma\cap\varLambda)$
\,($\gamma\in\varGamma_X$), hence, by Proposition
\ref{pr:Gibbs|cond}, we can rewrite (\ref{eq:Psi3}) in the form
\begin{equation}\label{eq:Psi4}
\gcl(A_\varLambda) =\int_{\varGamma_{\varLambda}} \mynn F(\xi)\mypp
S_\varLambda(\xi)\, \lambda_{\theta}(\rd\xi),
\end{equation}
with $S_\varLambda(\xi)\in L^1(\varGamma_\varLambda,
\lambda_\theta)$.
Therefore, in order to show that the right-hand side of
(\ref{eq:Psi4}) vanishes, it suffices to check that
\begin{equation}\label{eq:=0}
\int_{\varGamma_{\varLambda}} \mynn
F(\xi)\,\lambda_{\theta}(\rd\xi)=0.
\end{equation}

To this end, substituting here the definition  (\ref{eq:F}) and
changing the order of integration, we can rewrite the integral in
(\ref{eq:=0}) as
\begin{align*}
\int_{\frakX^2}
\theta^{\otimes\myp2}(B_\varLambda(\bar{y}_1,\bar{y}_2))
\,\eta(\rd\bar{y}_1)\,\eta(\rd\bar{y}_2),
\end{align*}
where
$$
B_\varLambda(\bar{y}_1,\bar{y}_2)):= \{(x_1,x_2)\in\varLambda^2:
\,x_1+y_1=x_2+y_2\ \text{\,for some }\, y_1\in\bar{y}_1,\
y_2\in\bar{y}_2\}.
$$
It remains to note that
\begin{align*}
\theta^{\otimes\myp2}\bigl(B_\varLambda(\bar{y}_1,\bar{y}_2)\bigr)
&=\int_\varLambda
\theta\left(\textstyle\bigcup\nolimits_{y_1\in\bar{y}_1}
\!\bigcup\nolimits_{y_2\in\bar{y}_2}
\{x_1+y_1-y_2\}\right)\theta(\rd x_1)\\[.2pc]
&\le \sum_{y_1\in\bar{y}_1}\int_\varLambda \theta
\left(\textstyle\bigcup\nolimits_{y_2\in\bar{y}_2}\mynn
\{x_1+y_1-y_2\}\right) \theta(\rd x_1)\\
&=\sum_{y_1\in\bar{y}_1}\int_\varLambda \theta
\bigl(D_{\{x_1+y_1\}}(\bar{y}_2)\bigr)\,\theta(\rd x_1) =0\qquad
(\eta\text{-a.s.}),
\end{align*}
since, by assumption (\ref{eq:condB2}),
$\theta\bigl(D_{\{x_1+y_1\}}(\bar{y}_2)\bigr)=0$ \,($\eta$-a.s.) and
$\bar{y}_1$ contains at most countably many coordinates. Hence,
(\ref{eq:=0}) follows and so part (b) is proved.
\endproof

\begin{remark}\label{rm:a}
In the Poisson cluster case (see \cite[Theorem 2.7\myp(a)]{BD3}),
conditions \mbox{(a-i)} and \mbox{(a-ii)} of Theorem
\ref{th:properClusterGibbs}\myp(a) are not only sufficient but also
necessary for the local finiteness of cluster configurations. While
it is obvious that condition \mbox{(a-i)} is always necessary, there
may be a question as to whether condition \mbox{(a-ii)} is such in
the case of a Gibbs cluster measure $\gcl$. Inspection of the proof
of Theorem \ref{th:properClusterGibbs}\myp(a) shows that the
difficulty here lies in the questionable relationship between the
conditions $\sigma(\mathcal{Z}_B)<\infty$ and
$\hatgamma(\mathcal{Z}_{B})<\infty$ ($\hatg$-a.s.) (which are
equivalent in the Poisson cluster case). According to Remark
\ref{rm:kappa<const} (see the Appendix), under the hypothesis of
boundedness of the first-order correlation function
$\kappa_{\myn\g}^1$, the former implies the latter, but the converse
may not always be true. Simple counter-examples can be constructed
by considering translation-invariant pair interaction potentials
$\varPhi(\{x_1,x_2\})=\phi_0(x_1-x_2)\equiv\phi_{0}(y-x)$ such that
$\phi_0(x)=+\infty$ on some subset $\varLambda_\infty\subset X$ with
$\theta(\varLambda_\infty)=\infty$. However, if $\kappa_{\myn\g}^1$
is bounded below and the mean number of configuration points in a
set $B$ is finite then the measure $\theta(B)$ must be finite (see
Remark~\ref{rm:kappa<const}).
\end{remark}

\begin{remark}\label{rm:b}
Similarly to Remark~\ref{rm:a}, it is of interest to ask whether
conditions \mbox{(b-i)} and \mbox{(b-ii)} of Theorem
\ref{th:properClusterGibbs}\myp(b) are necessary for the simplicity
of the cluster measure $\gcl$ (as is the case for the Poisson
cluster measure, see \cite[Theorem 2.7\myp(b)]{BD3}).
However, in the Gibbs cluster case this is not so; more precisely,
\mbox{(b-i)} is of course necessary, but \mbox{(b-ii)} may not be
satisfied. For a simple counter-example, let the in-cluster measure
$\eta$ be concentrated on a single-point configuration
$\bar{y}=(0)$, so that the droplet cluster $D_{\{x\}}(\bar{y})$ is
reduced to a single-point set $\{x\}$. Here, any measure $\theta$
with atoms will not satisfy condition \mbox{(b-ii)}. On the other
hand, consider a Gibbs measure $\g$ with a hard-core
translation-invariant pair interaction potential
$\varPhi(\{x_1,x_2\})=\phi_0(x_1-x_2)\equiv\phi_{0}(y-x)$, where
$\phi_0(x)=+\infty$ for $|x|<r_0$ and $\phi_0(x)=0$ for $|x|\ge
r_0$; then in each admissible configuration $\gamma$ any two points
are at least at a distance $r_0$, and in particular any such
$\gamma$ is simple.
\end{remark}

\begin{remark}\label{rm:b'}
As suggested by Remarks \ref{rm:a} and \ref{rm:b}, it is plausible
that conditions \mbox{(a-ii)} and \mbox{(b-ii)} of Theorem
\ref{th:properClusterGibbs} \textit{are} necessary for the claims
(a) and (b), respectively,
if the interaction potential of the underlying Gibbs measure $\g$ is
finite on all finite configurations, i.e., $\varPhi(\xi)<+\infty$
for all $\xi\in\varGamma_X^{\myp 0}$.
\end{remark}

In conclusion of this section, let us state some criteria sufficient
for conditions \mbox{(a-ii)} and \mbox{(b-ii)} of Theorem
\textup{\ref{th:properClusterGibbs}} (see details in
\cite[\S\myp2.4]{BD3}). Assume for simplicity that the in-cluster
configurations are a.s.-finite, $\eta\{N_X(\bar{y})<\infty\}=1$.

\begin{proposition}
\label{pr:a2} Either of the following conditions is sufficient for
\textup{(\ref{eq:condA2})}\textup{:}

\textup{(a-ii$'$)} \,For any compact set $B\in\calB(X)$, the
$\theta$-measure of\/ its translations is uniformly bounded,
\begin{equation}\label{sigma-cond}
C_B:=\sup_{x\in X} \theta(B+x)<\infty,
\end{equation}
and, moreover, the mean number of\/ in-cluster points is finite,
\begin{equation}\label{sigma-cond*}
\int_{\frakX} N_X(\bar{y})\,\eta(\rd\bar{y})<\infty.
\end{equation}

\textup{(a-ii$''$)} \,The coordinates of vector $\bar{y}$ are a.s.\
uniformly bounded, that is, there is a compact set
$B_0\in\mathcal{B}(X)$ such that \,$N_{X\setminus B_0}(\bar{y})=0$
for $\eta$-a.a.\ $\bar{y}\in\frakX$.
\end{proposition}

\begin{proposition}
\label{pr:b2} Either of the following conditions is sufficient for
\textup{(\ref{eq:condB2})}\textup{:}

\textup{(b-ii$'$)} \,The measure\/ $\theta$ is non-atomic, that is,
\,$\theta\{x\}=0$\myp{} for each\/ $x\in X$.

\textup{(b-ii$''$)} \,For each\/ $x\in X$, 
\,$N_{\{x\}}(\bar{y})=0$ \,for\, $\eta$-a.a.\ $\bar{y}\in\frakX$.
\end{proposition}

\section{Quasi-invariance and the integration-by-parts formula}
\label{sec:3}

From now on, we restrict ourselves to the case where
$X={\mathbb{R}}^{d}$. Henceforth, we assume that the in-cluster
configurations are a.s.-finite, $\eta\{N_X(\bar{y})<\infty\}=1$;
hence, the component $X^\infty$ representing infinite clusters (see
Section~\ref{sec:2.1}) may be dropped, so the set $\frakX$ is now
redefined as $\mathfrak{X}:=\bigsqcup_{\myp n\in\ZZ_+}\mynn X^{n}$
(cf.\ (\ref{eq:calX})). Note that condition \mbox{(a-i)} of Theorem
\ref{th:properClusterGibbs} is then automatically satisfied.

We assume throughout that the correlation function
$\kappa_{\myn\g}^1(x)$ is bounded, which implies by Theorem
\ref{th:properClusterGibbs} that the same is true for the
correlation function $\kappa_{\myn\hatg}^1(z)$. Let us also impose
conditions (\ref{sigma-cond}) and (\ref{sigma-cond*}) which, by
Proposition \ref{pr:a2}, ensure that condition \mbox{(a-ii)} of
Theorem \ref{th:properClusterGibbs}\myp(a) is fulfilled and so
$\gcl$-a.a.\ configurations $\gamma\in\varGamma_X^\sharp$ are
locally finite. According to Proposition \ref{prop1}, condition
(a-ii) also implies that $\sigma(\mathcal{Z}_{B})<\infty$ providing
that $\theta(B)<\infty$, where the set $\mathcal{Z}_{B}\subset
\mathcal{Z}$ is defined in (\ref{eq:K}).

Finally, we require the probability measure $\eta$ on $\frakX$ to be
absolutely continuous with respect to the Lebesgue measure
$\rd\bar{y}$,
\begin{equation}\label{eq:h}
\eta (\rd \bar{y})=h(\bar{y})\,{{\mathrm{d}}}\bar{y},\qquad \bar{y}
=(y_{1},\dots ,y_{n})\in X^{n}\quad (n\in\ZZ_+).
\end{equation}
By Proposition \ref{pr:b2}\myp\mbox{(b-ii$''$)}, this implies that
Gibbs CPP configurations $\gamma $
are $\gcl$-a.s.\ simple (i.e., have no multiple points). Altogether,
the above assumptions ensure that $\gcl$-a.a.\ configurations
$\gamma $ belong to the proper configuration space $\varGamma_{X}$.

Our aim in this section is to prove the quasi-invariance of the
measure $\gcl$ with respect to
compactly supported diffeomorphisms of $X$ (Section~\ref{sec:3.2}),
and to establish an integration-by-parts formula
(Section~\ref{sec:3.3}). We begin in Section~\ref{sec:3.1} with a
brief description of some convenient ``manifold-like'' concepts and
notations first introduced in \cite{AKR1} (see also
\cite[\S\myp4.1]{BD3}), which furnish a suitable framework for
analysis on configuration spaces.

\subsection{Differentiable functions on configuration spaces}
\label{sec:3.1}

Let $T_{x}X$ be the tangent space of $X={\mathbb{R}}^d$ at point
$x\in X$. It can be identified in the natural way with
${\mathbb{R}}^{d}$, with the corresponding (canonical) inner product
denoted by a ``fat'' dot~$\CD$\,. The gradient on $X$ is denoted by
$\nabla $. Following \cite{AKR1}, we define the ``tangent space'' of
the configuration space $\varGamma_{X}$ at $\gamma \in \varGamma
_{X}$ as the Hilbert space $T_{\gamma}\varGamma
_{X}:=L^{2}(X\rightarrow TX;\,{\rd}\gamma )$, or equivalently
$T_{\gamma}\varGamma_{X}=\bigoplus_{x\in \gamma}T_{x}X$. The scalar
product in $T_{\gamma}\varGamma_{X}$ is denoted by
$\langle\cdot,\cdot\rangle_{\gamma}$, with the corresponding norm
$|\cdot|_\gamma$. A vector field $V$ over $\varGamma_{X}$ is a map
$\varGamma_{X}\ni \gamma\mapsto V(\gamma )=(V(\gamma )_{x})_{x\in
\gamma}\in T_{\gamma}\varGamma_{X}$. Thus, for vector fields
$V_1,V_2$ over $\varGamma_{X}$ we have
\begin{equation*}
\left\langle V_1(\gamma ),V_2(\gamma)
\right\rangle_{\gamma}=\sum_{x\in \gamma}V_1(\gamma )_{x} \CD
V_2(\gamma)_{x},\qquad \gamma \in \varGamma_X.
\end{equation*}

For $\gamma\in\varGamma_{X}$ and $x\in\gamma$, denote by
${\mathcal{O}}_{\gamma,\myp x}$ an arbitrary open neighborhood of
$x$ in $X$ such that ${\mathcal{O}}_{\gamma ,\myp x}\cap \gamma
=\{x\}$. For any measurable function $F:\varGamma_{X}\rightarrow
{{\mathbb{R}}}$, define the function $F_{x}(\gamma,\cdot):
{\mathcal{O}}_{\gamma,\myp x}\to\RR$ by $F_{x}(\gamma,y):=F((\gamma
\setminus \{x\})\cup \{y\})$, and set
\begin{equation*}
\nabla_{\myn x} F(\gamma ):=\left.\nabla F_{x}(\gamma,y)
\right|_{y=x},\qquad x\in X,
\end{equation*}
provided that $F_{x}(\gamma,\cdot)$ is differentiable at $x$.

Recall that for a function $\phi:X\to\RR$ its support $\supp\phi$ is
defined as the closure of the set $\{x\in X\!: \phi(x)\ne 0\}$.
Denote by ${\mathcal{FC}}(\varGamma_{X})$ the class of functions on
$\varGamma_{X}$ of the form
\begin{equation}\label{local-funct}
F(\gamma)=f(\langle \phi_{1},\gamma \rangle,\dots,\langle
\phi_{k},\gamma\rangle),\qquad \gamma \in \varGamma_{X},
\end{equation}
where $k\in\NN$, \,$f\in C_{b}^{\infty}(\RR^{k})$ ($:=$ the set of
$C^{\infty}$-functions on ${\mathbb{R}}^{k}$ bounded together with
all their derivatives), and $\phi_{1},\dots ,\phi_{k}\in
C_{0}^{\infty}(X)$ ($:=$ the set of $C^{\infty }$-functions on $X$
with compact support). Each $F\in {\mathcal{FC}}(\varGamma_{X})$ is
local, that is, there is a compact $K\subset X$ (which may depend on
$F$) such that $ F(\gamma)=F(\gamma\cap K)$ for all
$\gamma\in\varGamma_{X}$. Thus, for a fixed $\gamma $ there are
finitely many non-zero derivatives $\nabla_{\myn x} F(\gamma)$.

For a function $F\in {\mathcal{FC}}(\varGamma_{X})$ its
$\varGamma$-gradient $\nabla^{\varGamma}\myn F\equiv
\nabla^{\varGamma}_{\! X} F$ is defined as
\begin{equation}\label{eq:G-gradient}
\nabla^{\varGamma}\myn F(\gamma ):=(\nabla_{\myn x} F(\gamma))_{x\in
\gamma}\in T_{\gamma}\varGamma_{X},\qquad \gamma \in \varGamma_{X},
\end{equation}
so the directional derivative of $F$ along a vector field $V$ is given by
\begin{equation*}
\nabla_{\myn V}^{\varGamma}F(\gamma ):=\langle
\nabla^{\varGamma}\myn F(\gamma ),V(\gamma
)\rangle_{\gamma}=\sum_{x\in \gamma }\nabla_{\myn x} F(\gamma) \CD
V(\gamma)_{x},\qquad \gamma \in \varGamma_{X}.
\end{equation*}
Note that the sum here contains only finitely many non-zero terms.

Further, let ${\mathcal{FV}}(\varGamma_{X})$ be the class of
cylinder vector fields $V$ on $\varGamma_{X}$ of the form
\begin{equation}\label{vf}
V(\gamma)_{x}=\sum_{i=1}^{k}A_{i}(\gamma)\mypp v_{i}(x)\in
T_{x}X,\qquad x\in X,
\end{equation}
where $A_{i}\in {\mathcal{FC}}(\varGamma_{X})$ and $v_{i}\in
\Vect_{0}(X)$ ($:=$ the space of compactly supported
$C^{\infty}$-smooth vector fields on $X$), \,$i=1,\dots ,k$ \,($k\in
\NN$). Any vector filed $v\in \Vect_{0}(X)$ generates a constant
vector field $V$ on $\varGamma_{X}$ defined by
$V(\gamma)_{x}:=v(x)$. We shall preserve the notation $v$ for it.
Thus,
\begin{equation}  \label{eq:grad-new}
\nabla_{\myn v}^{\varGamma} F(\gamma )=\sum_{x\in \gamma}
\nabla_{\myn x}F(\gamma)\CD v(x),\qquad \gamma \in \varGamma_X.
\end{equation}

The approach based on ``lifting'' the differential structure from
the underlying space $X$ to the configuration space $\varGamma_X$ as
described above can also be applied to the spaces
${\mathfrak{X}}=\bigsqcup_{\myp n=0}^{\infty }X^n$,
$\mathcal{Z}=X\times\mathfrak{X}$ and $\varGamma_{\mathfrak{X}}$,
$\varGamma_{\mathcal{Z}}$.
For these spaces, we will use the analogous notations as above
without further explanation.

\subsection{$\Diff_{0}$-quasi-invariance}\label{sec:3.2}

In this section, we discuss the property of quasi-invariance of the
measure $\gcl$ with respect to diffeomorphisms of $X$. Let us start
by describing how diffeomorphisms of $X$ act on configuration
spaces. For a measurable map $\varphi :X\to X$, its \textit{support}
$\supp\varphi$ is defined as the smallest closed set containing all
$x\in X$ such that $\varphi (x)\ne x$. Let ${\Diff}_{0}(X)$ be the
group of diffeomorphisms of $X$ with \textit{compact support}. For
any $\varphi \in {\Diff}_{0}(X)$, consider the corresponding
``diagonal'' diffeomorphism
$\bar{\varphi}:\mathfrak{X}\to\mathfrak{X}$ acting on each
constituent space $X^n$ ($n\in\mathbb{Z}_+$) as
\begin{equation}\label{eq:phi-hat0}
X^n\ni \bar{y}=(y_{1},\dots,y_{n})\mapsto
\bar{\varphi}(\bar{y}):=(\varphi (y_{1}),\dots ,\varphi(y_{n}))\in
X^n.
\end{equation}
For $x\in X$, we also define ``shifted'' diffeomorphisms
\begin{equation}\label{eq:phi-hat}
\bar{\varphi}_{x}(\bar{y}):=\bar\varphi(\bar{y}+x)-x,\qquad
\bar{y}\in {\mathfrak{X}}
\end{equation}
(see the shift notation (\ref{eq:shift})). Finally, we introduce a
special class of diffeomorphisms $\hat\varphi$ on $\mathcal{Z}$
acting only in the $\bar{y}$-coordinate as follows,
\begin{equation}\label{eq:hat-phi}
\hat\varphi(z):= (x,\bar{\varphi}_{x}(\bar{y}))\equiv
(x,\myp\bar\varphi(\bar{y}+x)-x),\qquad z=(x,\bar{y})\in
\mathcal{Z}.
\end{equation}

\begin{remark}\label{supp}
Note that, even though $K_\varphi:= \supp \varphi $ is compact in
$X$, the support of the diffeomorphism $\hat{\varphi}$ (again
defined as the closure of the set $\{z\in \mathcal{Z}:
\hat{\varphi}(z)\ne z\}$) is given by
$\supp\hat{\varphi}=\mathcal{Z}_{K_\varphi}$ (see~(\ref{eq:K})) and
hence is \textit{not} compact in the topology of $\mathcal{Z}$
(see Section~\ref{sec:2.1}).
\end{remark}

In the standard fashion, the maps $\varphi $ and $\hat{\varphi}$ can
be lifted to measurable ``diagonal'' transformations (denoted by the
same letters) of the configuration spaces $\varGamma_{X}$ and
$\varGamma_{\mathcal{Z}}$, respectively:
\begin{equation}\label{di*}
\begin{aligned} \varGamma_{X}\ni \gamma \mapsto \varphi (\gamma
):={}&\{\varphi (x),\ x\in \gamma \}\in \varGamma_{X},\\[.2pc]
\varGamma_{\mathcal{Z}}\ni \hatgamma\mapsto
\hat{\varphi}(\hatgamma):={}&\{\hat{\varphi}(z),\ (z)\in
\hatgamma\}\in \varGamma_{\mathcal{Z}}.
\end{aligned}
\end{equation}

The following lemma shows that the operator $\mathfrak{q}$ commutes
with the action of diffeomorphisms (\ref{di*}).\footnote{According
to relation (\ref{comm1}), $\mathfrak{q}$ is an \textit{intertwining
operator} between associated diffeomorphisms $\varphi$ and
$\hat{\varphi}$.}
\begin{lemma}
For any diffeomorphism $\varphi \in {\Diff}_0(X)$ and the
corresponding diffeomorphism $\hat{\varphi}$, it holds
\begin{equation}  \label{comm1}
\varphi\circ\mathfrak{q}=\mathfrak{q}\circ \hat{\varphi}.
\end{equation}
\end{lemma}

\proof The statement follows from the definition (\ref{eq:proj}) of
the map $\mathfrak{q}$ in view of the structure of diffeomorphisms
$\varphi$ and $\hat{\varphi}$ (see (\ref{eq:hat-phi}) and
(\ref{di*})).
\endproof

\begin{lemma}\label{lm:shift-free} The interaction potential\/
\,$\hat{\U}$\mypp{} defined in \textup{(\ref{eq:hatU})} is invariant
with respect to diffeomorphisms \textup{(\ref{eq:hat-phi})}, that
is, for any $\varphi\in\Diff_0(X)$ we have
$$
\hat{\U}(\hat{\varphi}(\hatgamma))=\hat{\U}(\hatgamma),
\qquad\hatgamma\in\varGamma_{\mathcal{Z}}.
$$
In particular, this implies the $\hat{\varphi}$-invariance of the
energy functionals defined in \textup{(\ref{eq:E})} and
\textup{(\ref{eq:E-E})}, that is, for any $\hat{\xi} \in
\varGamma_{\mathcal{Z}}^{\myp0}$ and $\hatgamma\in
\varGamma_{\mathcal{Z}}$,
\begin{equation*}
\hat{E}(\hat{\varphi}(\hat{\xi}\myp))=\hat{E}(\hat{\xi}\myp),\qquad
\hat{E}(\hat{\varphi}(\hat{\xi}\myp),\hat{\varphi}(\hatgamma))
=\hat{E}(\hat{\xi},\hatgamma).
\end{equation*}
\end{lemma} \proof
The claim readily follows by observing that a diffeomorphism
(\ref{eq:hat-phi}) acts on the $\bar{y}$-coordinates of points
$z=(x,\bar{y})$ in a configuration
$\hatgamma\in\varGamma_{\mathcal{Z}}$, while the interaction
potential $\hat{\U}$ (see (\ref{eq:hatU})) only depends on their
$x$-coordinates.
\endproof
As already mentioned (see (\ref{eq:h})), we assume that the measure
$\eta $ is absolutely continuous with respect to the Lebesgue
measure $\rd\bar{y}$
on ${\mathfrak{X}}$ and, moreover,
\begin{equation}\label{QI}
h(\bar{y}):=\frac{\eta (\rd\bar{y})}{\rd\bar{y}}>0\qquad {\text{for
\,a.a.}}\ \,\bar{y}\in {\mathfrak{X}}.
\end{equation}
This implies that the measure $\eta$ is quasi-invariant with respect
to the action of transformations
$\bar{\varphi}:{\mathfrak{X}}\rightarrow {\mathfrak{X}}$ \,($\varphi
\in {\Diff}_{0}(X)$), that is, for any $f\in\mathrm{M}_+(\frakX)$,
\begin{equation}\label{eq:eta-inv}
\int_{\frakX} f(\bar{y})\,\bar{\varphi}^*\eta(\rd\bar{y})
=\int_{\frakX} f(\bar{y})\,\rho_{\eta}^{\bar{\varphi}}(\bar{y})\,
\rd\bar{y},
\end{equation}
with the Radon--Nikodym density
\begin{equation}\label{density'}
\rho_{\eta}^{\bar{\varphi}}(\bar{y}):=\frac{\rd(\bar{\varphi}^*\eta)}{\rd\eta}(\bar{y})=
\frac{h(\bar{\varphi}^{-1}(\bar{y}))}{
h(\bar{y})}\,J_{\bar{\varphi}}(\bar{y})^{-1}
\end{equation}
(we set $\rho_{\eta}^{\bar{\varphi}}(\bar{y})=1$ if $h(\bar{y})=0$
or $h(\bar{\varphi}^{-1}(\bar{y}))=0$). Here
$J_{\bar{\varphi}}(\bar{y})$ is the Jacobian determinant of the
diffeomorphism $\bar{\varphi}$; due to the diagonal structure of
$\bar{\varphi}$ (see (\ref{eq:phi-hat0})) we have
$J_{\bar{\varphi}}(\bar{y})=\prod_{y_i\in\bar{y}}
J_{\varphi}(y_{i})$, where $J_{\varphi}(y)$ is the Jacobian
determinant of $\varphi$.

Due to the ``shift'' form of diffeomorphisms (\ref{eq:hat-phi}),
formulas (\ref{eq:eta-inv}), (\ref{density'}) readily imply that the
product measure
$\sigma(\rd{z})=\theta(\rd{x})\otimes\eta(\rd\bar{y})$ on
$\mathcal{Z}=X\times\frakX$ is quasi-invariant with respect to
$\hat\varphi$, that is, for each $\varphi\in\Diff_0(X)$ and any
$f\in\mathrm{M}_+(\mathcal{Z})$,
\begin{equation}\label{eq:sigma-inv}
\int_{\mathcal{Z}} f(z)\,\hat{\varphi}^*\sigma(\rd{z})
=\int_{\mathcal{Z}} f(z)\,\rho_{\varphi}(z)\, \sigma(\rd{z}),
\end{equation}
where the Radon--Nikodym density
$\rho_{\varphi}:=\rd(\hat\varphi^*\sigma)/\rd\sigma$ is given by
(see (\ref{density'}))
\begin{equation}\label{eq:rho-phi}
\rho_{\varphi}(z)= \rho_{\eta}^{\bar{\varphi}_x}(\bar{y})
\equiv\frac{h\bigl(\bar{\varphi}^{-1}(\bar{y}+x)-x\bigr)\,}{h(\bar{y})\,}
\,J_{\bar{\varphi}}(\bar{y}+x)^{-1},\quad
z=(x,\bar{y})\in\mathcal{Z}.
\end{equation}

We can now state our result on the quasi-invariance of the measure
$\hatg$.

\begin{theorem}\label{th:inv}
The Gibbs measure $\hatg$ constructed in Section
\textup{\ref{sec:2.4}} is quasi-invariant with respect to the action
of diagonal diffeomorphisms\/ $\hat{\varphi}$ on
$\varGamma_{\mathcal{Z}}$ \textup{(}$\varphi\in\Diff_0(X)$\textup{)}
defined by formula \textup{(\ref{eq:hat-phi})}, with the
Radon--Nikodym density
$R_{\hatg}^{\hat{\varphi}}=\rd(\hat{\varphi}^{\ast} \hatg)/\rd\hatg$
given by
\begin{equation}\label{RND}
R_{\hatg}^{\hat{\varphi}}(\hatgamma) = \prod_{z\in \hat{\gamma\myp}}
\myp\rho_{\varphi}(z),\qquad\hatgamma\in\varGamma_{\mathcal{Z}},
\end{equation}
where $\rho_{\varphi}(z)$ is defined in \textup{(\ref{eq:rho-phi})}.
\end{theorem}
\proof First of all, note that $\rho_{\varphi}(z)=1$ for any
$z=(x,\bar{y})\notin \supp\hat{\varphi}=\mathcal{Z}_{K_\varphi}$,
where $K_\varphi=\supp\varphi$ (see Remark \ref{supp}), and
$\sigma(\mathcal{Z}_{K_\varphi})<\infty$ by Proposition \ref{prop1}.
On the other hand, Theorem \ref{th:g-hat}\myp(b) implies that the
correlation function
$\kappa_{\myn\hatg}^1$
is bounded. Therefore, by Remark \ref{rm:kappa<const} (see the
Appendix) we obtain that $\hatgamma(\mathcal{Z}_{K_\varphi})<\infty$
for $\hatg$-a.a.\ configurations
$\hatgamma\in\varGamma_{\mathcal{Z}}$, hence the product in
(\ref{RND}) contains
finitely many terms different from $1$ and so the function
$R_{\hatg}^{\hat{\varphi}}(\hatgamma)$ is well defined. Moreover, it
satisfies the ``localization'' equality
\begin{equation}\label{eq:R=R}
R_{\hatg}^{\hat{\varphi}}(\hatgamma)=
R_{\hatg}^{\hat{\varphi}}(\hatgamma\cap
\mathcal{Z}_{K_\varphi})\qquad \text{for}\ \ \hatg\text{-a.a.}\
\,\hatgamma\in\varGamma_{\mathcal{Z}}.
\end{equation}

Following \cite[\S\myp2.8, Theorem 2.8.2]{KunaPhD},
the proof of the theorem will be based on the use of Ruelle's
equation (see the Appendix, Theorem \ref{th:NZR}). Namely, according
to (\ref{ruelle}) with $\varLambda= \mathcal{Z}_{K_\varphi}$, for
any function $F\in\mathrm{M}_+(\varGamma_{\mathcal{Z}})$ we have
\begin{align}
\notag &\int_{\varGamma_{\mathcal{Z}}}\mynn F(\hatgamma)
\,\hat{\varphi}^{\ast}{\hatg}(\rd\hatgamma)=\int_{\varGamma_{\mathcal{Z}}}
\mynn F(\hat{\varphi}(\hatgamma))
\,\hatg(\rd\hatgamma)\\
\notag &\ \
=\int_{\varGamma_{\varLambda}}\biggl(\int_{\varGamma_{\mathcal{Z}
\setminus \varLambda}} \! F(\hat{\varphi}(\hat{\xi} \cup
\hatgamma'))\,{\re}^{-\hat{E}(\hat{\varphi}(\hat{\xi}\myp))
-\hat{E}(\hat{\varphi}(\hat{\xi}\myp),\myp
\hat{\varphi}(\hatgamma'))}\,
\hatg(\rd\hatgamma')\biggr)\mypp\lambda_{\sigma}(\rd\hat{\xi}\myp) \\
\notag &\ \
=\int_{\varGamma_{\varLambda}}\biggl(\int_{\varGamma_{\mathcal{Z}\setminus
\varLambda}} \!F(\hat{\varphi}(\hat{\xi}\myp)\cup
\hatgamma')\,{\re}^{-\hat{E}(\hat{\varphi}(\hat{\xi}\myp))
-\hat{E}(\hat{\varphi}(\hat{\xi}\myp),\myp\hatgamma')}
\,\hatg(\rd\hatgamma')\biggr)\mypp\lambda_{\sigma}(\rd\hat{\xi}\myp)\\
\label{eq:bridge} &=\int_{\varGamma_{\varLambda}}\biggl(
\int_{\varGamma_{\mathcal{Z}\setminus \varLambda}} \!F(\hat{\xi}
\cup
\hatgamma')\,{\re}^{-\hat{E}(\hat{\xi}\myp)-\hat{E}(\hat{\xi},\myp
\hatgamma')}\,\hatg(\rd\hatgamma')\biggr)\,
\hat{\varphi}^{\ast}\lambda_{\sigma}(\rd\hat{\xi}\myp),
\end{align}
where $\lambda_\sigma$ is the Lebesgue--Poisson measure
corresponding to the reference measure $\sigma$ (see (\ref{eq:LP})).
Since $\sigma$ is quasi-invariant with respect to diffeomorphisms
$\hat\varphi$ (see (\ref{eq:sigma-inv})), it readily follows from
the definition (\ref{eq:LP}) that the restriction of the
Lebesgue--Poisson measure $\lambda_\sigma$ onto the set
$\varGamma_\varLambda$ is quasi-invariant with respect to
$\hat\varphi$, with the density given precisely by the function
(\ref{RND}). Hence, using the property (\ref{eq:R=R}), the
right-hand side of (\ref{eq:bridge}) is reduced to
\begin{align}
\notag
&\int_{\varGamma_{\varLambda}}\mynn\biggl(\int_{\varGamma_{\mathcal{Z}
\setminus \varLambda}}\!F(\hat{\xi} \cup
\hatgamma')\,{\re}^{-\hat{E}(\hat{\xi})-\hat{E}(\hat{\xi},\myp\hatgamma')}\,\hatg
(\rd\gamma^{\myp\prime})\mynn\biggr)\mypp R_{\hatg}^{\hat{\varphi}}
(\hat{\xi}\myp)\,\lambda_{\sigma}(\rd\hat{\xi}\myp) \\
\notag
&\qquad=\int_{\varGamma_{\varLambda}}\mynn\biggl(\int_{\varGamma_{\mathcal{Z}
\setminus \varLambda}}\!F(\hat{\xi} \cup
\hatgamma^{\myp\prime})\mypp R_{\hatg}^{\hat{\varphi}}(\hat{\xi}
\cup \hatgamma')\,{\re}^{-E(\hat{\xi})-E(\hat{\xi},
\hatgamma^{\myp\prime})}
\,{\hatg}(\rd\hatgamma^{\myp\prime})\mynn\biggr)\,\lambda_{\sigma}(\rd\hat{\xi}\myp) \\
\label{eq:R2} &\qquad\qquad=\int_{\varGamma_{\mathcal{Z}}}
\!F(\hatgamma)\mypp R_{\hatg}^{\hat{\varphi}}
(\hatgamma)\,{\hatg}(\rd\hatgamma),
\end{align}
where we have again used Ruelle's equation (\ref{ruelle}).

As a result, combining (\ref{eq:bridge}) and (\ref{eq:R2}) we obtain
\begin{equation}\label{eq:R3}
\int_{\varGamma_{\mathcal{Z}}}\mynn F(\hatgamma)
\,\hat{\varphi}^{\ast}{\hatg}(\rd\hatgamma)=\int_{\varGamma_{\mathcal{Z}}}
\!F(\hatgamma)\mypp R_{\hatg}^{\hat{\varphi}}
(\hatgamma)\,{\hatg}(\rd\hatgamma),
\end{equation}
which proves quasi-invariance of $\hatg$. In particular, setting
$F\equiv 1$ in (\ref{eq:R3}) yields $\int_{\varGamma_{\mathcal{Z}}}
\!R_{\hatg}^{\hat{\varphi}} (\hatgamma)\,{\hatg}(\rd\hatgamma)=1$,
and hence $R_{\hatg}^{\hat{\varphi}}\in
L^1(\varGamma_{\mathcal{Z}},\gcl)$.
\endproof

\begin{remark}
Note that the Radon--Nikodym density $R_{\hatg}^{\hat{\varphi}}$
defined by (\ref{RND}) does not depend on the background interaction
potential $\U$. As should be evident from the proof above, this is
due to the special ``shift'' form of the diffeomorphisms
$\hat{\varphi}$ (see (\ref{eq:hat-phi})) and the cylinder structure
of the interaction potential $\hat{\U}$ (see (\ref{eq:hatU})). In
particular, the expression (\ref{RND}) applies to the
``interaction-free'' case with $\U\equiv 0$ (and hence
$\hat{\U}\equiv0$), where the Gibbs measure
$\g\in\mathscr{G}(\theta,\U=0)$ is reduced to the Poisson measure
$\pi_\theta$ on $\varGamma_X$ with intensity measure $\theta$ (see
the Appendix), while the Gibbs measure
$\hatg\in\mathscr{G}(\sigma,\hat{\U}=0)$ amounts to the Poisson
measure $\pi_\sigma$ on $\varGamma_{\mathcal{Z}}$ with intensity
measure $\sigma$.
\end{remark}

\begin{remark}
As is essentially well known (see, e.g., \cite{AKR1,Sk}),
quasi-invariance of a Poisson measure on the configuration space
follows directly from the quasi-invariance of its intensity measure.
For a proof adapted to our slightly more general setting (where
diffeomorphisms are only assumed to have the support of finite
measure), we refer the reader to \cite[Proposition A.1]{BD3}.
Incidentally, the expression for the Radon--Nikodym derivative given
in \cite{BD3} (see also
\cite[Proposition 2.2]{AKR1}) contained a superfluous normalizing
constant, which in our context would read
\begin{equation*}
C_\varphi:=\exp\left(\int_{\mathcal{Z}}
\bigl(1-\rho_{\varphi}(z)\bigr)\,\sigma(\rd{z})\right)
\end{equation*}
(cf.\ (\ref{RND})). In fact, it is easy to see that $C_\varphi=1$;
indeed, $\rho_\varphi=1$ outside the set $\supp
\hat{\varphi}=\mathcal{Z}_{K_\varphi}$ with
$\sigma(\mathcal{Z}_{K_\varphi})<\infty$ (see Proposition
\ref{prop1}), hence
\begin{align*}
\ln C_\varphi=\int_{\mathcal{Z}_{K_\varphi}}
\bigl(1-\rho_{\varphi}(z)\bigr)\,\sigma(\rd{z})&=
\sigma(\mathcal{Z}_{K_\varphi})-\sigma(\hat{\varphi}^{-1}(\mathcal{Z}_{K_\varphi}))=0.
\end{align*}
\end{remark}

Let $\Q:L^{2}(\varGamma_{X},\gcl)\rightarrow
L^{2}(\varGamma_{\mathcal{Z}},\hatg)$ be the isometry defined by the
map ${\mathfrak{q}}$ (see (\ref{eq:proj})),
\begin{equation}\label{eq:I}
({\Q}F)(\hatgamma):=F\circ\mathfrak{q}(\hatgamma),\qquad \hatgamma
\in \varGamma_{\mathcal{Z}},
\end{equation}
and consider the corresponding adjoint operator
\begin{equation}\label{eq:I*}
\Q^*:\myp L^{2}(\varGamma_{\mathcal{Z}}, \hatg)\to
L^{2}(\varGamma_{X},\gcl).
\end{equation}

\begin{lemma}\label{lm:L1}
The operator $\mathcal{I}_{\mathfrak{q}}^{\ast}$ defined by
\textup{(\ref{eq:I*})} can be extended to the operator
\begin{equation*}
\Q^*:\myp L^{1}(\varGamma_{\mathcal{Z}},\hatg)\rightarrow
L^{1}(\varGamma _{X},\gcl).
\end{equation*}
\end{lemma}

\proof Note that $\Q$ can be viewed as a bounded operator acting
from $L^{\infty }(\varGamma _{X},\gcl)$ to
$L^{\infty}(\varGamma_{\mathcal{Z}},\hatg)$. This implies that the
adjoint operator $\Q^*$ is a bounded operator on the corresponding
dual spaces, $\Q^*:\myp
L^{\infty}(\varGamma_{\mathcal{Z}},\hatg)^{\prime}\to
L^{\infty}(\varGamma _{X},\gcl)^{\prime }$.

It is known (see \cite{GL}) that, for any sigma-finite measure space
$(M,\mu)$, the corresponding space $L^{1}(M,\mu)$ can be identified
with the subspace $V$ of the dual space $L^{\infty}(M,\mu)^{\prime}$
consisting of all linear functionals on $L^{\infty}(M,\mu)$
continuous with respect to bounded convergence in
$L^{\infty}(M,\mu)$. That is, $\ell\in V$ if and only if
$\ell(\psi_{n})\rightarrow 0$ for any $\psi_n\in L^{\infty }(M,\mu)$
such that $|\psi_{n}|\le 1$ and $\psi_{n}(x)\rightarrow 0$ as
$n\to\infty$ for $\mu$-a.a.\ $x\in M$. Hence, to prove the lemma it
suffices to show that, for any $F\in
L^{1}(\varGamma_{\mathcal{Z}},\hatg)$, the functional $\Q^*F\in
L^{\infty} (\varGamma_{\mathcal{Z}},\hatg)^{\prime}$ is continuous
with respect to bounded convergence in
$L^{\infty}(\varGamma_{\mathcal{Z}},\hatg)$. To this end, for any
sequence $(\psi_{n})$ in $L^{\infty}(\varGamma_{X},\gcl)$ such that
$|\psi_{n}|\le 1$ and $\psi_{n}(\gamma)\rightarrow 0$ for
$\gcl$-a.a.\ $\gamma\in \varGamma_{X}$, we have to prove that
$\Q^*F(\psi_{n})\to0$.

Let us first show that $\Q\psi_{n}(\hatgamma)\equiv
\psi_{n}(\mathfrak{q}(\hatgamma))\rightarrow 0$ for $\hatg$-a.a.\
$\hatgamma\in \varGamma_{\mathcal{Z}}$. Set
\begin{align*}
A_{\psi}:={}&\{\gamma\in\varGamma_{X}:\psi_{n}(\gamma)\rightarrow
0\}\in\mathcal{B}(\varGamma_{X}),\\[.3pc]
\hat{A}_\psi:={}&\{\hatgamma\in\varGamma_{\mathcal{Z}}:
\psi_{n}(\mathfrak{q}(\hatgamma))\rightarrow
0\}\in\mathcal{B}(\varGamma_{\mathcal{Z}}),
\end{align*}
and note that $\hat{A}_\psi=\mathfrak{q}^{-1}(A_\psi)$; then,
recalling the relation (\ref{eq:gcl*}),
we get
\begin{equation*}
\hatg(\hat{A}_\psi)=\hatg\bigl(\mathfrak{q}^{-1}(A_\psi)\bigr)=\gcl(A_\psi)=1,
\end{equation*}
as claimed. Now, by the dominated convergence theorem this implies
\begin{equation*}
\Q^*F(\psi_{n})=\int_{\varGamma_{\mathcal{Z}}}
F(\hatgamma)\,\Q\psi_{n} (\hatgamma)\,\hatg(\rd\hatgamma)\rightarrow
0,
\end{equation*}
and the proof is complete.
\endproof

Taking advantage of Theorem \ref{th:inv} and applying the projection
construction, we obtain our main result in this section.
\begin{theorem}\label{q-inv}
The Gibbs cluster measure $\gcl$ is quasi-invariant with respect to
the action of\/ ${\Diff}_{0}(X)$ on $\varGamma_{X}$. The
corresponding Radon--Nikodym density is given by
$R_{\gcl}^{\varphi}=\Q^* R_{\hatg}^{\hat{\varphi}}$.
\end{theorem}

\proof Note that, due to (\ref{eq:gcl*}) and (\ref{comm1}),
$$
\gcl\circ \varphi^{-1}=\hatg\circ
\mathfrak{q}^{-1}\circ\varphi^{-1}=\hatg\circ\hat{\varphi}^{-1}\circ
\mathfrak{q}^{-1}.
$$
That is, $\varphi^*\gcl=\gcl\circ\varphi^{-1}$ is a push-forward of
the measure $\hat{\varphi}^{\ast}\hatg=\hatg\circ\hat{\varphi}^{-1}$
under the map $\mathfrak{q}$, that is,
$\varphi^*\gcl=\mathfrak{q}^*\hat{\varphi}^{\ast}\hatg$. In
particular, if $\hat{\varphi}^{\ast}\hatg$ is absolutely continuous
with respect to $\hatg$ then so is $\varphi^{\ast}\gcl$ with respect
to $\gcl$. Moreover, by the change of measure (\ref{eq:gcl*}) and by
Theorem \ref{th:inv}, for any
$F\in L^\infty(\varGamma_X,
\gcl)$ we have
\begin{align}
\int_{\varGamma_{X}} \mynn
F(\gamma)\,\varphi^*\myn\gcl(\rd\gamma)
&=\int_{\varGamma_{\mathcal{Z}}}{\Q} F(\hatgamma)\,
\hat{\varphi}^{\ast}\hatg(\rd\hatgamma)
\label{eq:L1} 
=\int_{\varGamma_{\mathcal{Z}}}\mynn
{\Q}F(\hatgamma)\myp R_{\hatg}^{\hat{\varphi}}
(\hatgamma)\,\hatg(\rd\hatgamma).
\end{align}
By Lemma \ref{lm:L1}, the operator $\Q^*$ acts from
$L^1(\varGamma_{\mathcal{Z}},\hatg)$ to $L^1(\varGamma_X,\gcl)$.
Therefore, by the change of measure (\ref{eq:ghat}) the right-hand
side of (\ref{eq:L1}) can be rewritten as
$$
\int_{\varGamma_{X}} \mynn F(\gamma)\mypp
(\Q^*R_{\hatg}^{\hat{\varphi}})(\gamma )\,\gcl(\rd\gamma ),
$$
which completes the proof.
\endproof

\begin{remark}
The Gibbs cluster measure $\gcl$ on the configuration space
$\varGamma_{X}$ can be used to construct a unitary representation
$U$ of the diffeomorphism group ${\Diff}_{0}(X)$ by operators in
$L^{2}(\varGamma_{X},\gcl)$, given by the formula
\begin{equation}\label{eq:U}
U_{\varphi}F(\gamma )=\sqrt{R_{\gcl}^{\varphi}(\gamma)}\,F(\varphi
^{-1}(\gamma )),\qquad F\in L^{2}(\varGamma_{X},\gcl).
\end{equation}
Such representations, which can be defined for arbitrary
quasi-invariant measures on $\varGamma_{X}$, play a significant role
in the representation theory of the group ${\Diff}_{0}(X)$
\cite{I,VGG} and quantum field theory \cite{GGPS,Goldin}. An
important question is whether the representation (\ref{eq:U}) is
irreducible. According to \cite{VGG}, this is equivalent to the
${\Diff}_{0}(X)$-ergodicity of the measure $\gcl$, which in our case
is equivalent to the ergodicity of the measure $\hatg$ with respect
to the group of transformations $\hat{\varphi}$
\,($\varphi\in\Diff_{0}(X)$). Adapting the technique developed in
\cite{KS}, it can be shown that the aforementioned ergodicity of
$\hatg$ is valid if and only if
$\hatg\in\ext\mathscr{G}(\sigma,\hat{\U}\myp)$. In turn, the latter
is equivalent to $\g\in\ext\mathscr{G}(\theta,\U\myp)$, provided
that $\g\in\mathscr{G}_{\mathrm{R}}(\theta,\U\myp)$ (see
Corollary~\ref{cor:2.9}).
\end{remark}

\subsection{Integration-by-parts formula}\label{sec:3.3}

Let us first prove simple sufficient conditions for our measures on
configuration spaces to belong to the corresponding moment classes
$\mathcal{M}^n$ (see the Appendix, formula (\ref{eq:Mn})).

\begin{lemma}\label{lm:M^n}
\textup{(a)} \,Let\/ $\g\in \mathscr{G}(\theta,\U\myp)$, and suppose
that the correlation functions\/ $\kappa_{\myn\g}^{m}$ are
bounded for all\/ $m=1,\dots,n$. Then $\hatg\in
\mathcal{M}^{n}(\varGamma _{\mathcal{Z}})$, that is,
\begin{equation}\label{4.4}
\int_{\varGamma_{\mathcal{Z}}} |\langle
f,\hatgamma\rangle|^{n}\,\hatg(\rd\hatgamma)<\infty,\qquad f\in
C_0(\mathcal{Z}).
\end{equation}
Moreover, the bound \textup{(\ref{4.4})} is valid for any function
$f\in \bigcap_{\myp m=1}^{\myp n} \myn L^{m}(\mathcal{Z},\sigma)$.

\textup{(b)} \,If, in addition, the total number of components of a
random vector\/ $\bar{y}\in\frakX$ has a finite\/ $n$-th
moment\,\footnote{Cf.\ our standard assumption (\ref{sigma-cond*}),
where $n=1$.} with respect to the measure\/ $\eta$,
\begin{equation}\label{eq:N<}
\int_{\frakX} N_{X}(\bar{y})^{n}\,\eta(\rd\bar{y})<\infty,
\end{equation}
then\/ $\gcl\in \mathcal{M}^{n}(\varGamma _{X})$.
\end{lemma}

\proof (a) \,Using the multinomial expansion, for any $f\in
C_0(\mathcal{Z})$ we have
\begin{align}\notag
\int_{\varGamma_{\mathcal{Z}}} |\langle
f,\hatgamma\rangle|^{n}\,\hatg(\rd\hatgamma)
&\le\int_{\varGamma_{\mathcal{Z}}}
\left(\sum\nolimits_{z\in\hatgamma} |f(z)|\right)^{n}
\,\hatg(\rd\hatgamma)  \\
\label{eq:multi} &=\sum_{m=1}^{n}\int_{\varGamma_{\mathcal{Z}}}
\sum_{\{z_1\myn,\dots,\myp z_m\} \subset\hatgamma}
\phi_n(z_1,\dots,z_m)\,\hatg(\rd\hatgamma),
\end{align}
where $\phi_n(z_1,\dots,z_m)$ is a symmetric function given by
\begin{equation}\label{eq:varPsi}
\phi_n(z_1,\dots,z_m):=\sum_{\substack{ i_{1}\myn,\dots,\myp i_{m}\ge1 \\[.1pc]
 i_{1}+\dots+i_{m}=n}} \frac{n!}{i_1!\cdots i_m!}\,
|f(z_{1})|^{i_{1}}\mynn\cdots |f(z_{m})|^{i_{m}}.
\end{equation}
By the definition (\ref{corr-funct}),
the integral on the right-hand side of (\ref{eq:multi}) is reduced
to
\begin{equation}\label{eq:k=}
\frac{1}{m!}\int_{\mathcal{Z}^{m}}
\phi_n(z_{1},\dots,z_{m})\,\kappa_{\myn\hatg}^{m}(z_{1},\dots,z_{m})
\,\sigma (\rd{z}_{1})\cdots\sigma (\rd{z}_{m}).
\end{equation}
By Theorem \ref{th:g-hat}\myp(b), the hypotheses of the lemma imply
that $0\le \kappa_{\myn\hatg}^{m}\le a_{\hatg}$ ($m=1,\dots,n$) with
some constant $a_{\hatg}<\infty$. Hence, substituting
(\ref{eq:varPsi}) we obtain that the integral in (\ref{eq:k=}) is
bounded by
\begin{equation} \label{4.3}
a_{\hatg}\sum_{\substack{i_{1},\dots,\myp i_{m}\ge1 \\[.1pc]
i_{1}+\dots+\myp i_{m}=\myp n}} \frac{n!}{i_1!\cdots
i_m!}\prod_{j=1}^m \int_{\mathcal{Z}}|f(z_{j})|^{i_{j}}
\,\sigma(\rd{z}_{j})<\infty,
\end{equation}
since each integral in (\ref{4.3}) is finite owing to the assumption
$f\in C_0(\mathcal{Z})$. Moreover, the bound (\ref{4.3}) is valid
for any function $f\in \bigcap_{\myp m=1}^{\myp n}
L^{m}(\mathcal{Z},\sigma )$. Returning to (\ref{eq:multi}), this
yields (\ref{4.4}).


(b) \,Using the change of measure (\ref{eq:gcl*}), for any $\phi \in
C_{0}(X)$ we obtain
\begin{align}
\int_{\varGamma_{X}} |\langle\phi,\gamma
\rangle|^{n}\,\gcl(\rd\gamma) =\int_{\varGamma_{\mathcal{Z}}}
|\langle\phi,\mathfrak{q}(\hatgamma)\rangle|^{n}\,\hatg(\rd\hatgamma)
\label{eq:ff} =\int_{\varGamma_{\mathcal{Z}}}|\langle
\mathfrak{q}^*\phi,\hatgamma\rangle|^{n}\,\hatg(\rd\hatgamma),
\end{align}
where
\begin{equation}\label{eq:q*}
\mathfrak{q}^*\phi(x,\bar{y}):=\sum_{y_i\in\bar{y}}\phi
(y_{i}+x),\qquad (x,\bar{y})\in\mathcal{Z}.
\end{equation}
Due to part (a) of the lemma, it suffices to show that
$\mathfrak{q}^{\ast }\phi \in L^{m}(\mathcal{Z},\sigma)$ for any
$m=1,\dots,n$. By the elementary inequality $(a_1+\cdots+a_k)^m\le
k^{m-1}(a_1^m+\cdots+a_k^m)$, from (\ref{eq:q*}) we have
\begin{equation}\label{eq:q*1}
\int_{\mathcal{Z}}|\mathfrak{q}^{\ast}\phi (z)|^{m}\,\sigma (\rd{z})
\leq \int_{\mathcal{Z}} N_X(\bar{y})^{m-1}\sum_{y_{i}\in \bar{y}}
|\phi (y_{i}+x)|^{m}\,\sigma(\rd{x}\times\rd\bar{y}).
\end{equation}
Recalling that $\sigma=\theta\otimes\eta$ and denoting
$b_\phi:=\sup_{x\in X}|\phi(x)|<\infty$ and
$K_\phi:=\supp\phi\subset X$, the right-hand side of (\ref{eq:q*1})
is dominated by
\begin{align*}
\int_{\frakX} N_X(\bar{y})^{m-1}&\left(\myn (b_\phi)^m\sum_{y_{i}\in
\bar{y}} \int_X \mathbf{1}_{K_\phi-y_i}(x)\,\theta(\rd{x})\right)\eta (\rd\bar{y})\\
&=(b_\phi)^m \int_{\frakX} N_X(\bar{y})^{m-1}\sum_{y_{i}\in
\bar{y}} \theta(K_\phi-y_i)\;\eta (\rd\bar{y})\\
&\le (b_\phi)^m\,\sup_{y\in X}\theta (K_\phi-y) \int_{\frakX}
N_X(\bar{y})^{m}\,\eta (\rd\bar{y})<\infty,
\end{align*}
according to the assumptions (\ref{sigma-cond}) and (\ref{eq:N<}).
\endproof

In the rest of this section, we shall assume that the conditions of
Lemma \ref{lm:M^n} are satisfied with $n=1$. Thus, the measures
$\g$, $\hatg$ and $\gcl$ belong to the corresponding
$\mathcal{M}^{1}$-classes.

Let $v\in \Vect_{0}(X)$ (:= the space of compactly supported smooth
vector fields on $X$), and define a vector field $\hat{v}_{x}$ on
${\mathfrak{X}}$ by the formula
\begin{equation}\label{eq:v-hat}
\hat{v}_{x}(\bar{y}):=\left( v(y_{1}+x),\dots ,v(y_{n}+x)\right)
,\qquad \bar{y}=(y_{1},\dots ,y_{n})\in {\mathfrak{X}}.
\end{equation}
Observe that the measure $\eta $ satisfies the following
integration-by-parts formula,
\begin{equation}\label{eq:etaIBP}
\int_{{\mathfrak{X}}}\nabla^{\hat{v}_{x}} \mynn
f(\bar{y})\,\eta(\rd\bar{y})=-\int_{\frakX}
f(\bar{y})\,\beta^{\hat{v}}_{\eta}(x,\bar{y})\,\eta
(\rd\bar{y}),\qquad f\in C_{0}^{\infty}({\mathfrak{X}}),
\end{equation}
where $\nabla^{\hat{v}_{x}}$ is the derivative along the vector
field $\hat{v}_{x}$ and
\begin{equation}\label{eq:log-der}
\beta^{\hat{v}}_{\eta}(x,\bar{y}):=(\beta_{\eta}(\bar{y}),
\hat{v}_{x}(\bar{y}))_{T_{\bar{y}}{\mathfrak{X}}}+\Div\hat{v}_{x}(\bar{y})
\end{equation}
is the logarithmic derivative of
$\eta(\rd\bar{y})=h(\bar{y})\,\rd\bar{y}$ along $\hat{v}_{x}$,
expressed in terms of the vector logarithmic derivative
\begin{equation}\label{eq:nabla-h}
\beta_{\eta}(\bar{y}):=\frac{\nabla h(\bar{y})}{h(\bar{y})},\qquad
\bar{y} \in {\mathfrak{X}}.
\end{equation}

Let us define the space $H^{1\myn,\myp n}(\mathfrak{X)}$
($n\ge 1$) as the set of functions $f\in L^{n}(\frakX,\rd\bar{y})$
satisfying the condition
\begin{equation}\label{eq:Sobolev}
\int_{\mathfrak{X}}\left(\myn\sum_{y_i\in\bar{y}}
\,\bigl|\nabla_{\myn
y_i}f(\bar{y})\bigr|\right)^{\!n}\rd\bar{y}<\infty.
\end{equation}
Note that $H^{1\myn,\myp n}(\mathfrak{X)}$ is a linear space, due to
the elementary inequality $(|a|+|b|)^n\le
2^{n-1}\bigl(|a|^n+|b|^n\bigr)$.
\begin{lemma}\label{lm:B^n}
Assume that\/ $h^{1/n}\in H^{1\myn,\myp n}(\mathfrak{X)}$ for some
integer\/ $n\ge1$, and let the condition\/ \textup{(\ref{eq:N<})}
hold. Then\/ $\beta_{\eta}^{\hat{v}}\in L^{m}(\mathcal{Z},\sigma)$
for any\/ $m=1,\dots,n$.
\end{lemma}

\proof Firth of all, observe that the condition $h^{1/n}\in
H^{1\myn,\myp n}(\mathfrak{X)}$ implies that $h^{1/m}\in
H^{1\myn,\myp m}(\mathfrak{X)}$ for any $m=1,\dots,n$. Indeed,
$h^{1/m}\in L^{m}({\mathfrak{X},\rd\bar{y}})$ if and only if $h\in
L^{1}({\mathfrak{X},\rd\bar{y}})$; furthermore, using the definition
(\ref{eq:nabla-h}) of $\beta_{\eta}(\bar{y})$ we see that
\begin{align}
\notag \int_{\mathfrak{X}}\left(\myn\sum_{y_i\in\bar{y}}\,
\bigl|\nabla_{\myn y_i}
\myn\bigl(h(\bar{y})^{1/m}\bigr)\bigr|\right)^{\!m}\rd\bar{y}&=m^{-m}
\int_{\mathfrak{X}}\left(\myn\sum_{y_i\in\bar{y}}
\frac{|\nabla_{\myn y_i}
h(\bar{y})|}{h(\bar{y})^{1-1/m}}\right)^{\!m}\rd\bar{y}\\
\label{eq:h-eta} &= m^{-m}
\int_{\mathfrak{X}}\left(\myn\sum_{y_i\in\bar{y}}
|\beta_{\eta}(\bar{y})_i|\right)^{\!m}\eta(\rd\bar{y})<\infty,
\end{align}
since $\eta$ is a probability measure and hence
$L^{m}(\frakX,\eta)\subset L^{n}(\frakX,\eta)$ ($m=1,\dots,n$).

To show that $\beta_{\eta}^{\hat{v}}\in L^{m}(\mathcal{Z},\sigma)$,
it suffices to check that each of the two terms on the right-hand
side of (\ref{eq:log-der}) belongs to $L^{m}(\mathcal{Z},\sigma)$.
Denote $b_v:=\sup_{x\in X}\myn|v(x)|<\infty$, $K_v:=\supp v\subset
X$, and recall that $C_{K_v}:=\sup_{y\in X} \theta(K_v-y)<\infty$ by
condition (\ref{sigma-cond}). Using (\ref{eq:v-hat}), we have
\begin{align}
\notag \int_{\mathcal{Z}}
&|(\beta_{\eta}(\bar{y}),\hat{v}_{x}(\bar{y}))|^{m}
\,\sigma(\rd{x}\times\rd\bar{y}) \le\int_{\mathcal{Z}}\! \left(\myn
\sum_{y_i\in\bar{y}}|\beta_{\eta}(\bar{y})_{i}|
\cdot|v(y_{i}+x)|\right)^{\!m}\myn\theta(\rd{x})\,\eta(\rd\bar{y}) \\
\notag &\leq
(b_v)^{m-1}\int_{\mathfrak{X}}\left(\myn\sum_{y_i\in\bar{y}}
|\beta_{\eta}(\bar{y})_{i}|\right)^{\!m-1}\sum_{y_i\in\bar{y}}
|\beta_{\eta}(\bar{y})_{i}|\left(\int_{X}
|v(y_{i}+x)|\,\theta(\rd{x})\right)\eta(\rd\bar{y}) \\
\notag &\leq
(b_v)^{m}\int_{\mathfrak{X}}\left(\myn\sum_{y_i\in\bar{y}}
|\beta_{\eta}(\bar{y})_{i}|\right)^{\!m-1}\sum_{y_i\in\bar{y}}
|\beta_{\eta}(\bar{y})_{i}|\,\theta(K_v-y_i)\,\eta(\rd\bar{y}) \\
\label{eq:beta<infty} &\leq (b_v)^{m}\,C_{K_v}
\int_{\mathfrak{X}}\left(\myn\sum_{y_i\in\bar{y}}
|\beta_{\eta}(\bar{y})_{i}|\right)^{\!m} \eta(\rd\bar{y}) <\infty,
\end{align}
according to (\ref{eq:h-eta}). Similarly, denoting $d_v:=\sup_{x\in
X}|\Div v(x)|<\infty$ and again using (\ref{eq:v-hat}), we obtain
\begin{align}
\notag \int_{\mathcal{Z}}\bigl|\Div {}&
\hat{v}_{x}(\bar{y})\bigr|^{m}\,\sigma(\rd{x}\times\rd\bar{y})
=\int_{\mathcal{Z}}\left(\myn \sum_{y_i\in\bar{y}}\bigl|(\Div
v)(y_{i}+x)
\bigr|\right)^{\!m}\theta(\rd{x})\,\eta(\rd\bar{y}) \\
\notag &\leq (d_v)^{m-1}\int_{\mathfrak{X}}
N_X(\bar{y})^{m-1}\left(\myn\sum_{y_i\in\bar{y}}
\int_{X}\bigl|(\Div v)(y_{i}+x)\bigr|\,\theta(\rd{x})\right)\eta (\rd\bar{y}) \\
\notag &\leq (d_v)^{m}\int_{\mathfrak{X}}
N_X(\bar{y})^{m-1}\sum_{y_i\in\bar{y}}
\theta(K_v-y_i)\,\eta (\rd\bar{y}) \\
\label{eq:div<infty} &\leq (d_v)^{m}\,C_{K_v}\int_{\mathfrak{X}}
N_X(\bar{y})^{m}\,\eta(\rd\bar{y})<\infty,
\end{align}
according to the assumption (\ref{eq:N<}). As a result, combining
the bounds (\ref{eq:beta<infty})
and (\ref{eq:div<infty}), we see that $\beta_{\eta}^{\hat{v}}\in
L^{m}(\mathcal{Z},\sigma)$, as claimed.
\endproof


The next two theorems are our main results in this section.

\begin{theorem}\label{IBP-}
For any function $F\in {\mathcal{FC}}(\varGamma_{X})$, the Gibbs
cluster measure $\gcl$ satisfies the integration-by-parts formula
\begin{equation}\label{IBP0-}
\int_{\varGamma_{X}}\sum_{x\in \gamma}\nabla_{\myn x}F(\gamma)\CD
v(x)\,\gcl(\rd\gamma )=-\int_{\varGamma_{X}}F(\gamma )\,
B_{\myn\gcl}^{v}(\gamma )\,\gcl(\rd\gamma ),
\end{equation}
where $B_{\myn\gcl}^{v}(\gamma ):=\Q^{*}\langle
\beta^{\hat{v}}_{\eta}, \hatgamma\rangle\in
L^{1}(\varGamma_{X},\gcl)$ and \strut{}$\beta^{\hat{v}}_{\eta}$ is
the logarithmic derivative defined in \textup{(\ref{eq:log-der})}.
\end{theorem}

\proof For any function $F\in {\mathcal{FC}}(\varGamma_{X})$ and
vector field $v\in \Vect_{0}(X)$, let us denote for brevity
\begin{equation}\label{eq:H}
H(x,\gamma ):=\nabla_{\myn x}F(\gamma)\CD v(x), \qquad x\in X, \ \
\gamma\in\varGamma_X.
\end{equation}
Furthermore, setting $\hat{F}=\Q F:\varGamma_{\mathcal{Z}}\to\RR$ we
introduce the notation
\begin{equation}\label{eq:H-hat}
\hat{H}(z,\hatgamma):= \nabla_{\myn\bar{y}}\hat{F}(\hatgamma) \CD
\hat{v}_{x}(\bar{y}),\qquad z=(x,\bar{y})\in\mathcal{Z},\ \
\hatgamma\in\varGamma_{\mathcal{Z}}.
\end{equation}
From these definitions, it is clear that
\begin{equation}\label{eq:QHH}
{\Q}\Biggl(\mynn\sum_{x\in \gamma}
H(x,\gamma)\Biggr)(\hatgamma)=\sum_{z\in
\hatgamma}\hat{H}(z,\hatgamma),\qquad
\hatgamma\in\varGamma_{\mathcal{Z}}.
\end{equation}

Let us show that the Gibbs measure $\hatg$ on
$\varGamma_{\mathcal{Z}}$ satisfies the following
integration-by-parts formula:
\begin{equation}\label{IBPp}
\int_{\varGamma_{\mathcal{Z}}}\sum_{z\in
\hat{\gamma\myp}}\hat{H}(z,\hatgamma)\,\hatg(\rd\hatgamma)
=-\int_{\varGamma_{\mathcal{Z}}}
\!\hat{F}(\hatgamma)B_{\myn\hatg}^{\hat{v}}(\hatgamma)\,
\hatg(\rd\hatgamma),
\end{equation}
where the logarithmic derivative
$B_{\myn\hatg}^{\hat{v}}(\hatgamma):=
\langle\beta^{\hat{v}}_{\eta},\hatgamma\rangle$ belongs to
$L^{1}(\varGamma_{\mathcal{Z}},\hatg)$ (by Lemmas
\ref{lm:M^n}\myp(b) and \ref{lm:B^n} with $n=1$). By the change of
measure (\ref{eq:ghat}) and due to relation (\ref{eq:QHH}), we have
\begin{align}
\notag \int_{\varGamma_{\mathcal{Z}}}\myn\sum_{z\in \hatgamma\myp}
|\hat{H}(z,\hatgamma)|\;\hatg(\rd\hatgamma) &=\int_{\varGamma
_{X}}\sum_{x\in \gamma}|H(x,\gamma)|\;
\gcl(\rd\gamma ) \\
\label{eq:|H|} &\le \sup_{(x,\myp\gamma)}
|H(x,\gamma)|\int_{\varGamma _{X}}\sum_{x\in\gamma}
\mathbf{1}_{K_{v}}(x)\:\gcl(\rd\gamma),
\end{align}
where $K_v:=\supp v$ is a compact set in $X$. Note that the
right-hand side of (\ref{eq:|H|}) is finite, since the function $H$
is bounded (see (\ref{eq:H})) and, by Lemma \ref{lm:M^n}\myp(b),
$\gcl\in\mathcal{M}^{1}(\varGamma_{X})$. Therefore, by Remark
\ref{rm:NZ} we can apply Nguyen--Zessin's equation (\ref{eq:NZ})
with the function $\hat{H}$ to obtain
\begin{equation}\label{eq:HH}
\int_{\varGamma_{\mathcal{Z}}}\myn\sum_{z\in \hat{\gamma\myp}}
\hat{H}(z,\hatgamma)\,\hatg(\rd\hatgamma)
=\int_{\varGamma_{\mathcal{Z}}}\!\left(\int_{\mathcal{Z}}
\hat{H}(z,\hatgamma\cup\{z\})\,\re^{-\hat{E}(\{z\},\myp\hat{\gamma\myp})}\,
\sigma (\rd{z})\right)\hatg(\rd\hatgamma).
\end{equation}
Inserting the definition (\ref{eq:H}), using Lemma
\ref{lm:shift-free} and recalling that $\sigma=\theta\otimes\eta$
(see (\ref{eq:sigma})), let us apply the integration-by-parts
formula (\ref{eq:etaIBP}) for the measure $\eta$ to rewrite the
internal integral in (\ref{eq:HH}) as
\begin{align*}
&\int_{X}{\re}^{-E(\{x\},\mypp p_{X}(\hat{\gamma\myp}))}\left(
\int_{{\mathfrak{X}}}\nabla _{\bar{y}}
\hat{F}(\hatgamma\cup\{(x,\bar{y})\})\CD \hat{v}_{x}(\bar{y})\,\eta
(\rd\bar{y})\right)\theta(\rd x)\\
&\qquad=-\int_{X} \re^{-E(\{x\},\myp p_{X}(\hat{\gamma\myp}))}
\left(\int_{{\mathfrak{X}}}\hat{F}(\hatgamma\cup\{(x,\bar{y})\})\mypp
\beta^{\hat{v}}_{\eta}(x,\bar{y})\, \eta (\rd\bar{y})\right) \theta
(\rd x)\\
&\qquad\qquad=-\int_{\mathcal{Z}} \re^{-E(\{p_{X}(z)\},\mypp
p_{X}(\hat{\gamma\myp}))} \hat{F}(\hatgamma\cup\{z)\})\mypp
\beta^{\hat{v}}_{\eta}(z)\, \sigma(\rd{z}).
\end{align*}
Returning to (\ref{eq:HH}) and again using Nguyen--Zessin's equation
(\ref{eq:NZ}), we see that the right-hand side of (\ref{eq:HH}) is
reduced to
\begin{align*}
-\int_{\varGamma_{\mathcal{Z}}} \mynn\sum_{z\in
\hat{\gamma\myp}}\hat{F}(\hatgamma)\mypp\beta^{\hat{v}}_{\eta}(z)\,
\hatg(\rd\hatgamma)&=-\int_{\varGamma_{\mathcal{Z}}} \mynn
\hat{F}(\hatgamma)\mypp B_{\myn\hatg}^{\hat{v}}\,\hatgamma\rangle\,
\hatg(\rd\hatgamma),
\end{align*}
which proves formula (\ref{IBPp}).

Now, using equality (\ref{eq:QHH}), we obtain
\begin{align*}
\int_{\varGamma_{X}}\sum_{x\in \gamma}H(x,\gamma)\,\gcl (\rd\gamma )
&=\int_{\varGamma_{\mathcal{Z}}}\!\Biggl(\mynn\sum_{(x,\myp\bar{y})\in
\hat {\gamma}}\!\nabla_{\myn\bar {y}}{\Q}F(\hatgamma)\CD
\hat{v}_{x}(\bar {y})\mynn\Biggr)\,
\hatg(\rd\hatgamma)\\
&=-\int_{\varGamma_{\mathcal{Z}}}{\Q}F(\hatgamma)\mypp
B_{\myn\hatg}^{\hat{v}}
(\hatgamma)\;\hatg(\rd\hatgamma) \\
&=-\int_{\varGamma_{X}}F(\gamma)\,
\Q^*B_{\myn\hatg}^{\hat{v}}(\gamma)\,\gcl(\rd\gamma),
\end{align*}
where $\Q^*B_{\myn\hatg}^{\hat{v}}\in L^1(\varGamma_X,\gcl)$ by
Lemma~\ref{lm:L1}. Thus, formula (\ref{IBP0-}) is proved.
\endproof

\begin{remark}
Observe that the logarithmic derivative
$B_{\myn\hatg}^{\hat{v}}=\langle\beta^{\hat{v}}_{\eta},\hatgamma\rangle$
(see (\ref{IBPp})) does not depend on the interaction potential
\mypp$\U$, and in particular coincides with that in the case
$\U\equiv0$, where the Gibbs measure $\g$ is reduced to the Poisson
measure $\pi_\theta$. Nevertheless, the logarithmic derivative
$B_{\myn\gcl}^v$ does depend on $\U$ via the map $\Q^*$.
\end{remark}

\begin{remark}
Note that in Theorem \ref{IBP-} the reference measure $\theta$
does not have to be differentiable with respect to $v$.
\end{remark}

According to Theorem \ref{IBP-}, $B_{\gcl}^{v}\in
L^1(\varGamma_{\mathcal{Z}},\gcl)$. However, under the conditions of
Lemma \ref{lm:B^n} with $n\ge2$, this statement can be enhanced.

\begin{lemma}\label{lm:B^n*}
Assume that\/ $h^{1/n}\in H^{1\myn,\myp n}(\mathfrak{X)}$ for some
integer\/ $n\ge2$, and let the condition\/ \textup{(\ref{eq:N<})}
hold. Then\/ $B_{\myn\gcl}^{v}\mynn\in
L^{n}(\varGamma_{\mathcal{Z}},\gcl)$.
\end{lemma}

\proof By Lemmas \ref{lm:M^n}\myp(a) and \ref{lm:B^n}, it follows
that $\langle \beta^{\hat{v}}_{\eta}, \hatgamma\rangle\in
L^{n}(\varGamma_{\mathcal{Z}},\hatg)$. Let $s:=n/(n-1)$, so that
$n^{-1}+s^{-1}=1$. Note that $\Q$ can be treated as a bounded
operator acting from $L^{s}(\varGamma_{X},\gcl)$ to
$L^{s}(\varGamma_{\mathcal{Z}},\hatg)$. Hence, $\Q^{\ast}$ is a
bounded operator from
$L^{s}(\varGamma_{\mathcal{Z}},\hatg)^{\prime}=L^{n}(\varGamma_{\mathcal{Z}},\hatg)$
to $L^{s}(\varGamma_{X},\gcl)^{\prime}=L^{n}(\varGamma _{X},\gcl)$,
which implies that $B_{\myn\gcl}^{v}=\Q^*\langle
\beta^{\hat{v}}_{\eta}, \hatgamma\rangle\in
L^{n}(\varGamma_{\mathcal{Z}},\gcl)$.
\endproof

Formula (\ref{IBP0-}) can be extended to more general vector fields
on $ \varGamma_{X}$. Let ${\mathcal{FV}}(\varGamma_{X})$ be the
class of vector fields $V$ of the form $V(\gamma )=(V(\gamma
)_{x})_{x\in \gamma}$,
\begin{equation*}
V(\gamma )_{x}=\sum_{j=1}^{N}G_{j}(\gamma)\,v_{j}(x)\in T_{x}X,
\end{equation*}
where $G_{j}\in {\mathcal{FC}}(\varGamma_{X})$ and $v_{j}\in
\Vect_{0}(X)$, \,$j=1,\dots ,N$. For any such $V$ we set
\begin{equation*}
B_{\myn\gcl}^{V}(\gamma):=(\Q^{\ast}B_{\myn\hatg}^{{\Q}V})(\gamma),
\end{equation*}
where $B_{\myn\hatg}^{{\Q}V}(\hatgamma)$ is the logarithmic
derivative of $\hatg$ along ${\Q}V(\hat{
\gamma}):=V({\mathfrak{q}}(\hatgamma))$ (see \cite{AKR1}). Note that
${\Q}V$ is a vector field on $\varGamma_{\mathcal{Z}}$ owing to the
obvious equality
\begin{equation*}
T_{\hat{\gamma\myp}}\varGamma_{\mathcal{Z}}
=T_{{\mathfrak{q}}(\hat{\gamma\myp})}\varGamma_{X}.
\end{equation*}
Clearly,
\begin{equation*}
B_{\myn\gcl}^{V}(\gamma )=\sum_{j=1}^{N}\biggl(
G_{j}(\gamma)B_{\myn\gcl}^{v_{j}}(\gamma )+\sum_{x\in
\gamma}\nabla_{\myn x} G_{j}(\gamma )\CD v_{j}(x)\biggr) .
\end{equation*}

\begin{theorem}\label{IBP1}
For any $F_{1},F_{2}\in {\mathcal{FC}}(\varGamma_{X})$
and\/ $V\in {\mathcal{FV}}(\varGamma_{X})$, we have
\begin{equation*}
\begin{aligned}
\int_{\varGamma_{X}}&\sum_{x\in \gamma}\nabla_{\myn
x}F_{1}(\gamma)\CD
V(\gamma)_{x}\,F_{2}(\gamma )\;\gcl(\rd\gamma) \\
& =-\int_{\varGamma_{X}}F_{1}(\gamma)\,\sum_{x\in
\gamma}\nabla_{x}F_{2}(\gamma )\CD V(\gamma)_{x}\ \gcl
(\rd\gamma)\\
&\qquad -\int_{\varGamma_{X}}F_{1}(\gamma)
F_{2}(\gamma)B_{\myn\gcl}^{V}(\gamma )\,\gcl(\rd\gamma ).
\end{aligned}
\end{equation*}
\end{theorem}

\proof The proof can be obtained by a straightforward generalization
of the arguments used in the proof of Theorem \ref{IBP-}. \endproof

We define the \emph{vector logarithmic derivative} of $\gcl$ as a
linear operator
\begin{equation*}
B_{\myn\gcl}\!:\,{\mathcal{FV}}(\varGamma_{X})\rightarrow
L^{1}(\varGamma _{X},\gcl)
\end{equation*}
via the formula
\begin{equation*}
B_{\myn\gcl}\myn V(\gamma):=B_{\myn\gcl}^{V}\myn(\gamma).
\end{equation*}
This notation will be used in the next section.

\section{The Dirichlet form and equilibrium stochastic dynamics}\label{sec:4}

Throughout this section, we assume that the conditions of Lemma
\ref{lm:M^n} are satisfied with $n=2$. Thus, the measures $\g$,
$\hatg$ and $\gcl$ belong to the corresponding
$\mathcal{M}^{2}$-classes. Our considerations will involve the
$\varGamma $-gradients (see Section~\ref{sec:3.1}) on different
configuration spaces, such as $\varGamma _{X}$, $\varGamma
_{\mathfrak{X}}$ and $\varGamma_{\mathcal{Z}}$; to avoid confusion,
we shall denote them by $\nabla_{\mynn X}^{\varGamma }$,
$\nabla_{\mathfrak{X}}^{\varGamma}$ and
$\nabla_{\mathcal{Z}}^{\varGamma}$, respectively.

\subsection{The Dirichlet form associated with $\protect\gcl$}
\label{sec:4.1}

Let us introduce a pre-Dirichlet form ${\mathcal{E}}_{\myn\gcl}$
associated with the Gibbs cluster measure $\gcl$, defined on
functions $F_1,F_2\in\mathcal{FC}(\varGamma_{X})\subset
L^{2}(\varGamma_{X},\gcl)$ by
\begin{equation}\label{eq:E-mu}
{\mathcal{E}}_{\myn\gcl}(F_1,F_2):=\int_{\varGamma_{X}}\langle
\nabla_{\!X}^{\varGamma} F_1(\gamma),\nabla_{\!X}^{\varGamma}
F_2(\gamma)\rangle_{\gamma} \,\gcl(\rd\gamma).
\end{equation}
Let us also consider the operator $H_{\myn\gcl}$ defined
by
\begin{equation}\label{eq:Hcl}
H_{\myn\gcl}F:=-\Delta
^{\varGamma}F+B_{\myn\gcl}\!\nabla_{\!X}^{\varGamma} F,\qquad F\in
{\mathcal{FC}}(\varGamma_{X}),
\end{equation}
where $\Delta ^{\varGamma}F(\gamma ):=\sum_{x\in
\gamma}\Delta_{x}F(\gamma)$.

The next theorem readily follows from the general theory of
(pre-)Dirichlet forms associated with measures from the class
$\mathcal{M}^{2}(\varGamma_{X})$ (see \cite{AKR2,MR}).

\begin{theorem}
\textup{(a)} \,The pre-Dirichlet form \textup{(\ref{eq:E-mu})} is
well defined, i.e., ${\mathcal{E}}_{\myn\gcl}(F_{1},F_{2})<\infty $
for all $F_{1},F_{2}\in {\mathcal{FC}}(\varGamma_{X})$\textup{;}

\textup{(b)} \,The expression \textup{(\ref{eq:Hcl})} defines a
symmetric operator $H_{\myn\gcl}\!$ in $L^{2}(\varGamma_{X},\gcl)$
whose domain includes ${\mathcal{FC}}(\varGamma_{X})$\textup{;}

\textup{(c)} \,The operator\/ $H_{\myn\gcl}\mynn$ is the generator
of the pre-Dirichlet form\/ ${\mathcal{E}}_{\myn\gcl}$, i.e.,
\begin{equation}\label{generator}
{\mathcal{E}}_{\myn\gcl}(F_{1},F_{2})=\int_{\varGamma_{X}}
F_{1}(\gamma)\,H_{\myn\gcl}\myn F_{2}(\gamma)\,\gcl(\rd\gamma),\
\quad F_{1},F_{2}\in {\mathcal{FC}}(\varGamma_{X}).
\end{equation}
\end{theorem}

Formula (\ref{generator}) implies that the form
${\mathcal{E}}_{\myn\gcl}$ is closable. It follows from the
properties of the \textit{carr\'e du champ} \,$\sum_{x\in \gamma}
\!\nabla_{\myn{}x} F_{1}(\gamma)\CD \nabla_{\myn x}F_{2}(\gamma )$
that the closure of ${\mathcal{E}}_{\myn\gcl}$ (for which we shall
keep the same notation) is a quasi-regular local Dirichlet form on a
bigger state space $\overset{\,..}{\varGamma}_{X}$ consisting of all
integer-valued Radon measures on $X$ (see \cite{MR}). By the general
theory of Dirichlet forms (see \cite{MR0}), this implies the
following result (cf.\ \cite{AKR1,AKR2,BD3}).

\begin{theorem}\label{th:7.2}
There exists a conservative diffusion process
$\mathbf{X}=(\mathbf{X}_t,\,t\ge0)$ on
$\overset{\,\myp..}{\varGamma}_{X}$, properly associated with the
Dirichlet form $\mathcal{E}_{\myn\gcl}$, that is, for any function
$F\in L^{2}(\overset{\,\myp..}{\varGamma}_{ X},\gcl)$ and all\/
$t\ge0$, the map
\begin{equation*}
\overset{\,\myp..}{\varGamma}_{X}\ni \gamma \mapsto p_{t}F(\gamma)
:=\int_{\varOmega} F(\mathbf{X}_{t})\,\rd P_{\gamma}
\end{equation*}
is an\/ $\mathcal{E}_{\myn\gcl}$-quasi-continuous version of\/
$\exp(-tH_{\myn\gcl}) F$. Here $\varOmega$ is the canonical sample
space \textup{(}of $\overset{\,\myp..}{\varGamma}_X$-valued
continuous functions on $\mathbb{R}_+$\textup{)} and
$(P_\gamma,\,\gamma\in\overset{\,\myp..}{\varGamma}_X)$ is the
family of probability distributions of the process $\mathbf{X}$
conditioned on the initial value $\gamma=\mathbf{X}_0$. The process
$\mathbf{X}$ is unique up to $\gcl$-equivalence. In particular,
$\mathbf{X}$ is $\gcl$-symmetric \textup{(}i.e., $\int F_1\mypp
p_{t}F_2\,\rd\gcl = \int F_2\, p_{t} F_1\,\rd\gcl$ for all
measurable functions $F_1,F_2:\overset{\,\myp..}{\varGamma}_{
X}\to\mathbb{R}_{+}$\textup{)} and $\gcl$ is its invariant measure.
\end{theorem}

\subsection{Irreducibility of the Dirichlet form}

Similarly to (\ref{eq:E-mu}), let ${\mathcal{E}}_{\myn\hatg}$ be the
pre-Dirichlet form associated with the Gibbs measure $\hatg$,
defined on functions
$F_1,F_2\in\mathcal{FC}(\varGamma_{\mathcal{Z}})\subset L^{2}(
\varGamma_{\mathcal{Z}},\hatg)$ by
\begin{equation}\label{generator1}
{\mathcal{E}}_{\myn\hatg}(F_1,F_2):=\int_{\varGamma_{\mathcal{Z}}}\langle
\nabla_{\!\mathcal{Z}}^{\varGamma}\myn
F_1(\hatgamma),\nabla_{\!\mathcal{Z}}^{\varGamma}\myn
F_2(\hatgamma)\rangle_{\hat{\gamma\myp}} \,\hatg(\rd\hatgamma).
\end{equation}
The integral on the right-hand side of (\ref{generator1}) is well
defined because $\hatg\in\mathcal{M}^{2}\subset \mathcal{M}^{1}$.
The latter fact also implies that the gradient operator
$\nabla_{\myn\mathcal{Z}}^{\varGamma}$ can be considered as an
(unbounded) operator $L^{2}(\varGamma_{\mathcal{Z}},\hatg)\to
L^{2}V(\varGamma_{\mathcal{Z}},\hatg)$ with domain
$\mathcal{FC}(\varGamma_{\mathcal{Z}})$, where
$L^{2}V(\varGamma_{\mathcal{Z}},\hatg)$ is the space of
square-integrable vector fields on $\varGamma_{\mathcal{Z}}$. Since
the form ${\mathcal{E}}_{\myn\hatg}$ belongs to the class
$\mathcal{M}^2$, it is closable \cite{AKR2} (we keep the same
notation for the closure and denote by
$\mathcal{D}(\mathcal{E}_{\myn\hatg})$ its domain).

Our aim is to study a relationship between the forms
$\mathcal{E}_{\myn\gcl}$ and $\mathcal{E}_{\myn\hatg}$ and to
characterize in this way the kernel of $\mathcal{E}_{\myn\gcl}$. We
need some preparations. Let us recall that the projection map
$\mathfrak{q}:\mathcal{Z}\to\varGamma^\sharp_X$ was defined in
(\ref{proj1}) as $\mathfrak{q}:=\mathfrak{p}\circ s$, where
\begin{equation*}
s:\mathcal{Z}\ni (x,\bar{y})\mapsto \bar{y}+x\in \mathfrak{X}.
\end{equation*}
As usual, we preserve the same notations for the induced maps of the
corresponding configuration spaces. It follows directly from the
definition (\ref{eq:pr}) of the map $\mathfrak{p}$ that
\begin{equation}\label{commut0}
(\nabla_{\mynn X}^{\varGamma }F) \circ \mathfrak{p}=\nabla
_{\mathfrak{X}}^{\varGamma}(F\circ \mathfrak{p}),\qquad
F\in\mathcal{FC}(\varGamma_{X}),
\end{equation}
where we use the identification of the tangent spaces
\begin{equation}\label{tangent1}
T_{\bar{\gamma}}\varGamma_{\mathfrak{X}}
={\textstyle\bigoplus\limits_{\bar{y}\in\mathfrak{\bar{\gamma}}}}
\,T_{\bar{y}}\mathfrak{X}={\textstyle{\bigoplus\limits_{\bar{y}\in
\mathfrak{\bar{\gamma}}}}}\,{\textstyle\bigoplus\limits_{y_i\in\bar{y}}}\,T_{y_i}X
={\textstyle\bigoplus\limits_{y_i\in\mathfrak{p}(\bar{\gamma})}}
T_{y_i} X=T_{\mathfrak{p}(\bar{\gamma})} X.
\end{equation}

\begin{theorem}\label{th:4.3}
For the Dirichlet forms\/ $\mathcal{E}_{\myn\gcl}$ and\/
$\mathcal{E}_{\myn\hatg}$ defined in \textup{(\ref{generator})} and
\textup{(\ref{generator1})}, respectively, their domains satisfy the
relation\/ $\myp\Q(\mathcal{D}(\mathcal{E}_{\myn\gcl}))\subset
\mathcal{D}(\mathcal{E}_{\myn\hatg})$. Furthermore, $F\in \Ker
{\mathcal{{\mathcal{E}}}}_{\myn\gcl}$ if and only if\/
$\mathcal{I}_{\mathfrak{q}} F\in \Ker \mathcal{E}_{\myn\hatg}$\myp.
\end{theorem}

\proof Let us introduce a map
$\rd{s}^{\ast}\mynn:\,\mathfrak{X}\rightarrow \mathcal{Z}$ by the
formula
\begin{equation*}
\rd{s}^{\ast}(\bar{y}):=\left(\textstyle{\sum_{y_{i}\in \bar{y}}
y_{i}},\,\bar{y}\right),\qquad \bar{y}\in\frakX.
\end{equation*}
As suggested by the notation, this map coincides with the adjoint of
the derivative
\begin{equation*}
\rd{s}(z)\mynn:\,T_{z}\mathcal{Z}\rightarrow T_{s(z)}\mathfrak{X}
\end{equation*}
under the identification $T_{\bar{y}}\mathfrak{X}=\mathfrak{X}$ and
$T_{z}\mathcal{Z}=\mathcal{Z}$. A direct calculation shows that for
any differentiable function $f$ on $\mathfrak{X}$ the following
commutation relation holds:
\begin{equation}\label{commut1}
(\rd{s}^{\ast }\nabla f)\circ s=\nabla (f\circ s).
\end{equation}
Here the symbol $\nabla $ denotes the gradient on the corresponding
space (i.e., $\mathfrak{X}$ on the left and $\mathcal{Z}$ on the
right).

Let
\begin{equation*}
\rd{s}^{\ast}(\hatgamma):\,T_{s(\hatgamma)}\varGamma_{\frakX}
=\textstyle{\bigoplus\limits_{\bar{y}\in
s(\hatgamma)}}T_{\bar{y}}\mathfrak{X} \rightarrow
\textstyle{\bigoplus\limits_{z\in \hatgamma}}\,T_{z}\mathcal{Z}
=T_{\hatgamma}\varGamma_{\mathcal{Z}}
\end{equation*}
be the natural lifting of the operator $\rd{s}^{\ast}$. Further,
using (\ref{tangent1}), it can be interpreted as the operator
\begin{equation*}
\rd{s}^{\ast}(\hatgamma):\,T_{\mathfrak{q}(\hatgamma)}
\varGamma_{X}\rightarrow T_{\hatgamma}\varGamma_{\mathcal{Z}},
\end{equation*}
which induces the (bounded) operator
\begin{equation}\label{oper}
\Q\,\rd{s}^{\ast}:\,L^{2}V(\varGamma _{X},\gcl)\rightarrow
L^{2}V(\varGamma_{\mathcal{Z}},\hatg)
\end{equation}
acting according to the formula
\begin{equation*}
(\Q\,\rd{s}^{\ast}V)
(\hatgamma)=\rd{s}^{\ast}(\hatgamma)V(\mathfrak{q}(\hatgamma)),\qquad
V\in L^{2}V(\varGamma _{X},\gcl).
\end{equation*}
Formula (\ref{commut1}) together with (\ref{commut0}) implies that
\begin{equation*}
\left(\rd{s}^*\nabla _{\!X}^{\varGamma}F\right) \circ
\mathfrak{q}=\nabla _{\!\mathcal{Z}}^{\varGamma}(F\circ
\mathfrak{q}),\qquad F\in \mathcal{FC}(\varGamma_{X}),
\end{equation*}
or, in terms of operators acting in the corresponding
$L^{2}$-spaces,
\begin{equation*}
\Q\,\rd{s}^*\nabla_{\!X}^{\varGamma
}F=\nabla_{\!\mathcal{Z}}^{\varGamma}\Q F,\qquad F\in
\mathcal{FC}(\varGamma_{X}).
\end{equation*}
Therefore, for any $F\in \mathcal{FC}(\varGamma_{X})$
\begin{align}
\mathcal{E}_{\myn\hatg}(\Q F,\Q F) &=\int_{\varGamma_{\mathcal{Z}}}
|(\Q\,\rd{s}^{\ast
}\nabla_{\!X}^{\varGamma}F)(\hatgamma)|_\hatgamma^{2}
\:\hatg(\rd\hatgamma)  \label{comrel1} \\
\notag &=\int_{\varGamma_{\mathcal{Z}}}|\rd{s}^{\ast}
\nabla_{\!X}^{\varGamma}F(\mathfrak{q}(\hatgamma))|_\hatgamma ^{2}\:\hatg(\rd\hatgamma)\\
\notag &=\int_{\varGamma_{X}}|\rd{s}^*\nabla_{\!X}^{\varGamma}
F(\gamma)|_\gamma^{2}\:\gcl(\rd\gamma)\\
\label{sqcl} &\geq \int_{\varGamma_{X}} |\nabla_{\!X}^{\varGamma}
F(\gamma)|_\gamma^{2}\:\gcl(\rd\gamma)= \mathcal{E}_{\myn\gcl}(F,F),
\end{align}
where in (\ref{sqcl}) we used the obvious inequality
$|\rd{s}^{\ast}(\bar{y})| \geq |\bar{y}|$ \,($\bar{y}\in
\mathfrak{X}$).
Hence,
\begin{align*}
\Vert F\Vert_{{\mathcal{E}}_{\gcl}}^{2}&
:={\mathcal{E}}_{\myn\gcl}(F,F)+\int_{\varGamma_{X}}F^{2}\,{\rd\gcl} \\
& \leq {\mathcal{E}}_{\myn\hatg}(\Q F,\Q
F)+\int_{\varGamma_{{\mathcal{Z}}}}(\Q F)^{2}\,\rd\hatg=\Vert \Q
F\Vert_{{\mathcal{E}}_{\myn\hatg}}^{2},
\end{align*}
which implies that $\Q(\mathcal{D}(\mathcal{E}_{\myn\gcl}))\subset
\mathcal{D}({\mathcal{E}}_{\myn\hatg})$, thus proving the first part
of the theorem.

Further, using approximation arguments and continuity of the
operator (\ref{oper}), one can show that the equality
(\ref{comrel1}) extends to the domain
$\mathcal{D}({\mathcal{E}}_{\myn\gcl})$,
\begin{equation}\label{eq:QQ}
\mathcal{E}_{\myn\hatg}(\Q F,\Q F)= \int_{\varGamma_{\mathcal{Z}}}
|\Q\,\rd{s}^*
\nabla_{\!X}^{\varGamma}F|_{\hatgamma}^{2}\:\rd\hatg,\qquad
F\in\mathcal{D}(\mathcal{E}_{\myn\gcl}).
\end{equation}
Since $\Ker(\Q\,\rd{s}^{\ast}) =\{0\}$, formula (\ref{eq:QQ})
readily implies that $\Q F\in \Ker \mathcal{E}_{\myn\hatg}$ if and
only if $\nabla_{\!X}^{\varGamma}F=0$. In turn, due to equality
(\ref{sqcl}), the latter is equivalent to $F\in
\Ker\mathcal{E}_{\myn\gcl}$.
\endproof

Let us recall that a Dirichlet form $\mathcal{E}$ is called
\textit{irreducible} if the condition $\mathcal{E}(F,F)=0$ implies
that $F=\const$.

\begin{corollary} The Dirichlet form\/
${\mathcal{E}}_{\myn\gcl}$ is irreducible if\/
${\mathcal{E}}_{\myn\hatg}$ is so.
\end{corollary}

\proof Follows immediately from Theorem \ref{th:4.3} and the obvious
fact that if ${\mathcal{I}_{\mathfrak{q}}}F=\const$ ($\hatg$-a.s.)
then $F=\const$ ($\gcl$-a.s.).
\endproof

\begin{remark}
It follows from the general theory of Gibbs measures (see,
e.g., \cite{AKR2}) that the form $\mathcal{E}_{\myn\hatg}$ is
irreducible if and only if
$\hatg\in\ext\mathscr{G}(\sigma,\hat{\U}\myp)$, which is in turn
equivalent to $\g\in\ext\mathscr{G}(\theta ,\U\myp)$ (provided that
$\g\in\mathscr{G}_{\mathrm{R}}(\theta,\U\myp)$, see
Corollary~\ref{cor:2.9}).
\end{remark}


\section*{Acknowledgments}\label{sec:Ack}

The authors would like to thank Sergio Albeverio, Yuri Kondratiev,
Eugene Lytvynov and Tobias Kuna for
helpful discussions.

\appendix

\section{Gibbs measures on configuration spaces}

\renewcommand{\thesection}{A}



Let us briefly recall the definition and some properties of (grand
canonical) Gibbs measures on the configuration space
$\varGamma_{X}$. For a more systematic exposition and further
details, see the classical books \cite{Ge79,Preston,Ruelle}; more
recent useful references include \cite{AKR2,KunaPhD,KKS98}.

Denote by $\varGamma_{X}^{\myp0}:=
\{\gamma\in\varGamma_X:\,\gamma(X)<\infty\}$
the subspace of finite configurations in $X$. Let $\varPhi:
\varGamma_X^{\myp 0}\to \RR\cup \{+\infty \}$ be a measurable
function (called the \emph{interaction potential}) such that
$\varPhi(\emptyset)=0$.
A simple, most common example is that of the \textit{pair
interaction potential}, i.e., such that $\varPhi(\gamma)=0$ unless
configuration $\gamma$ consists of exactly two points.

\begin{definition}
The \emph{energy} $E:\varGamma_{X}^{\myp0}\rightarrow \mathbb{R}\cup
\{+\infty\}$ is defined by
\begin{equation}\label{eq:E}
E(\xi):=\sum_{\zeta\subset \xi}\varPhi(\zeta)\qquad (\xi \in
\varGamma_{X}^{\myp0}),\qquad E(\emptyset ):=0.
\end{equation}
The \emph{interaction energy} between configurations \myp$\xi \in
\varGamma_{X}^{\myp0}$ and $\gamma \in \varGamma_{X}$ is given by
\begin{equation}\label{eq:E-E}
E(\xi ,\gamma ):=\left\{
\begin{array}{ll}
\displaystyle\sum_{\gamma\supset\gamma'\in\varGamma_{X}^{\myp0}}
\!\varPhi(\xi\cup\gamma') & \displaystyle\quad \text{if}\
\sum_{\gamma\supset\gamma'\in\varGamma_{X}^{\myp0}}
\!|\myp\varPhi(\xi\cup\gamma')| <\infty , \\[1.7pc]
\displaystyle\ \ +\infty  & \quad \text{otherwise}.
\end{array}
\right.
\end{equation}
\end{definition}

\begin{definition}
Let $\mathscr{G}(\theta,\U\myp)$ denote the class of all
\textit{grand canonical Gibbs measures} corresponding to the
reference measure $\theta$ and the interaction potential $\U$, that
is, the probability measures on $ \varGamma_{X}$ that satisfy the
\textit{Dobrushin--Lanford--Ruelle \textup{(}DLR\textup{)} equation}
(see, e.g., \cite[Eq.~(2.17), p.~251]{AKR2}).
\end{definition}

In the present paper, we use an equivalent characterization of Gibbs
measures based on the following theorem, first proved by Nguyen and
Zessin \cite[Theorem~2]{NZ}.

\begin{theorem}\label{th:NZR}
A measure $\g$ on the configuration space $\varGamma_X$ belongs to
the Gibbs class $\mathscr{G}(\theta,\U\myp)$ if and only if either
of the following conditions holds\textup{:}

\textup{(i)} \textup{(}\emph{Nguyen--Zessin's equation}\textup{)}
\,For any function $H\in\mathrm{M}_+(X\times\varGamma_X)$,
\begin{equation}\label{eq:NZ}
\!\int_{\varGamma_{X}}\sum_{x_i\in\gamma} H(x_i,\gamma
)\,\g(\rd\gamma)=\int_{\varGamma_{X}}\!\biggl(\int_{X}H(x,\gamma
\cup \{x\})\,\re^{-E(\{x\},\myp\gamma )}\,\theta
(\rd{x})\mynn\biggr)\,\g(\rd\gamma).
\end{equation}

\textup{(ii)} \textup{(}\emph{Ruelle's equation}\textup{)} \,For any
bounded function $F\in\mathrm{M}_+(\varGamma_{X})$ and any compact
set $\varLambda\subset X$,
\begin{equation}\label{ruelle}
\int_{\varGamma_{X}}\! F(\gamma)\,\g(\rd\gamma)=
\int_{\varGamma_{\varLambda}}\re^{-E(\xi)}\biggl(\int_{\varGamma_{X\setminus
\varLambda}} \!F(\xi \cup \gamma')\
\re^{-E(\xi,\myp\gamma')}\,\g(\rd\gamma')\mynn\biggr)\,
\lambda_{\theta}(\rd\xi),
\end{equation}
where $\lambda_{\theta}$ is the \emph{Lebesgue--Poisson measure} on
$\varGamma_{X}^{\myp0}$ defined
by the formula
\begin{equation}\label{eq:LP}
\lambda_{\theta}({\rd}\xi)=\sum_{n=0}^{\infty}
\mathbf{1}\{\xi(\varLambda)=n\}\,\frac{1}{n!}
\,{\textstyle\bigotimes\limits_{x_i\in\xi}}\,\theta(\rd{x}_{i}),
\qquad \xi\in\varGamma^{\myp 0}_\varLambda.
\end{equation}
\end{theorem}

\begin{remark}\label{rm:NZ}
Using a standard argument based on the decomposition
$H=H^{+}-H^{-}$, \,$|H|=H^{+}+H^{-}$ with $H^{+}\myn:=\max\{H,0\}$,
$H^{-}\myn:=\max\{-H,0\}$, one can see that
equation (\ref{eq:NZ}) is also valid for an arbitrary measurable
function $H:X\times\varGamma_X\to\RR$ provided that
\begin{equation}\label{abs-int}
\int_{\varGamma_{X}}\sum_{x_i\in\gamma}|H(x_i,\gamma)|\:\g(\rd\gamma)<\infty.
\end{equation}
\end{remark}

\begin{remark}
In the original paper \cite{NZ}, the authors proved the result of
Theorem \ref{th:NZR} under additional assumptions of
\textit{stability} of the interaction potential $\U$ and
\textit{temperedness} of the measure $\g$. In subsequent work by
Kuna
\cite[Theorems 2.2.4, A.1.1]{KunaPhD}, these assumptions have been
removed.
\end{remark}

\begin{remark}\label{rm:apriori}
Inspection of \cite[Theorem~2]{NZ} or \cite[Theorem~A.1.1]{KunaPhD}
reveals that the proof of the implication $\textup{(\ref{eq:NZ})}
\Rightarrow \textup{(\ref{ruelle})}$ is valid for \textit{any} set
$\varLambda\in\mathcal{B}(X)$ satisfying a priori conditions
$\theta(\varLambda)<\infty$ and  $\gamma(\varLambda)<\infty$
($\g$-a.s.). Hence, Ruelle's equation (\ref{ruelle}) is valid for
such sets as well.
\end{remark}

In the ``interaction-free'' case where $\U\equiv 0$, the unique
grand canonical Gibbs measure coincides with the Poisson measure
$\pi_{\theta }$ (with intensity measure $\theta$). In the general
situation, there are various types of conditions to ensure that the
class $\mathscr{G}(\theta,\U\myp)$ is non-empty (see
\cite{Ge79,Preston,Ruelle} and also \cite{KunaPhD,K,KKS98}).

\begin{example} The
following are four classical examples of translation-invariant pair
interaction potentials (i.e., such that
$\U(\{x,y\})=\phi_{0}(x-y)\equiv\phi_{0}(y-x)$), for which
$\mathscr{G}(\theta,\U\myp)\ne\emptyset$.
\begin{enumerate}
\item[(1)]
(\textit{Hard core potential}) \,$\phi_0(x)=+\infty$ for $|x|\le
r_0$, otherwise $\phi_0(x)$ $=0$ \,($r_0>0$).
\item[(2)] (\textit{Purely repulsive potential})
\,$\phi_0\in C^2_0(\mathbb{R}^d)$, $\phi_0\ge0$ on $\mathbb{R}^d$,
and $\phi_0(0)>0$.
\item[(3)] (\textit{Lennard--Jones type potential}) \,$\phi_0\in
C^2(\mathbb{R}^d\setminus\{0\})$, $\phi_0\ge -a>-\infty$ on
$\mathbb{R}^d$, $\phi_0(x):=c|x|^{-\alpha}$ for $|x|\le r_1$
\,($c>0$, $\alpha>d$\myp), and \,$\phi_0(x)=0$ for $|x|>r_2$
\,($0<r_1<r_2<\infty$).
\item[(4)] (\textit{\emph{Lennard--Jones ``6--12'' potential}})
\,$d=3$, \,$\phi_0(x)=c(|x|^{-12}-|x|^{-6})$ for $x\ne0$ \,($c>0$)
and $\phi_0(0)=+\infty$.
\end{enumerate}
\end{example}

\begin{definition}
For a Gibbs measure $\g$ on $\varGamma_X$, its \textit{correlation
function $\kappa_{\myn\g}^n: X^n\to\RR_+$ of the $n$-th order}
($n\in\NN$) is defined by the following property: for any function
$\phi\in C_0(X^{n})$, symmetric with respect to permutations of its
arguments, it holds
\begin{multline}\label{corr-funct}
\int_{\varGamma_{X}}\sum_{\{x_{1}\myn,\dots,\myp x_{n}\} \subset
\gamma}
\phi(x_{1},\dots,x_{n})\,\g(\rd\gamma)\\
=\frac{1}{n!}\int_{X^{n}}\phi
(x_{1},\dots,x_{n})\,\kappa_{\myn\g}^{n}(x_{1},\dots,x_{n})\,\theta(\rd{x}_{1})\cdots\theta
(\rd{x}_{n}).
\end{multline}
\end{definition}
By a standard approximation argument, equation (\ref{corr-funct})
can be extended to any (symmetric) bounded measurable functions
$\phi:X^n\to\RR$ with support of finite $\theta^{\otimes
n}$-measure. For $n=1$ and $\phi(x)=\mathbf{1}_{B}(x)$,
the definition (\ref{corr-funct}) specializes to
\begin{equation}\label{corr-funct1}
\int_{\varGamma_{X}}\gamma(B)\,\g(\rd\gamma)=\int_{B}\kappa_{\myn\g}^{1}(x)\,\theta(\rd{x}).
\end{equation}
More generally, choosing $\phi(x_1,\dots,x_n)=\prod_{i=1}^n {\bf
1}_{B_i}(x_i)$ with arbitrary test sets $B_i\in\mathcal{B}(X)$, it
is easy to see that the definition (\ref{corr-funct}) is equivalent
to the following more explicit description (cf.\
\cite[p.~266]{AKR1}),
$$
\kappa_{\myn\g}^n(x_1,\dots,x_n)\,\theta(\rd{x}_1)\cdots\theta(\rd{x}_n)=
n!\cdot \g\{\gamma\in\varGamma_X: \gamma(\rd{x}_i)\ge 1,\
i=1,\dots,n\},
$$
also showing that indeed $\kappa_{\myn\g}^n\ge0$.

\begin{example}
In the Poisson case (i.e., $\U\equiv0$), we have
$\kappa_{\pi_\theta}^n(x)\equiv n!$ \,($n\in\NN$).
\end{example}


\begin{remark}\label{rm:NZ-kappa}
Using Nguyen--Zessin's equation (\ref{eq:NZ}) with
$H(x,\gamma)=\phi(x)$, from the definition (\ref{corr-funct1}) it
follows that
\begin{equation}\label{eq:corr1-NZ}
\kappa_{\myn\g}^1(x)=\int_{\varGamma_X}\re^{-E(\{x\},\gamma)}\,\g(\rd{\gamma}),\qquad
x\in X.
\end{equation}
In particular, the representation (\ref{eq:corr1-NZ}) implies that
if $\U\ge0$ (non-attractive interaction potential) then
$\kappa_{\myn\g}^1(x)\le 1$ for al $x\in X$, so that
$\kappa_{\myn\g}^1$ is bounded.
\end{remark}

\begin{remark}\label{rm:kappa<const}
If the first-order correlation function $\kappa_{\myn\g}^1(x)$ is
integrable on any set $B\in\mathcal{B}(X)$ of finite
$\theta$-measure (for instance, if $\kappa_{\myn\g}^1$ is bounded on
$X$, cf.\ Remark \ref{rm:NZ-kappa}), then, according to
(\ref{corr-funct1}), the mean number of points in $\gamma\cap B$ is
finite, also implying that $\gamma(B)<\infty$ for $\g$-a.a.\
configurations $\gamma\in\varGamma_X$ (cf.\
Remark~\ref{rm:apriori}). Conversely, if $\kappa_{\myn\g}^1$ is
bounded below (i.e., $\kappa_{\myn\g}^1(x)\ge c>0$ for all $x\in X$)
and the mean number of points in $\gamma\cap B$ is finite, then it
follows from (\ref{corr-funct1}) that $\theta(B)<\infty$.
\end{remark}

\begin{definition}\label{def:Mp}
For a probability measure $\mu$ on $\varGamma_{X}$, the notation
$\mu\in \mathcal{M}^{n}(\varGamma_{X})$ signifies that
\begin{equation}\label{eq:Mn}
\int_{\varGamma_{X}}|\langle \phi,\gamma\rangle|^{n}\,\mu
(\rd\gamma)<\infty,\qquad \phi \in C_{0}(X).
\end{equation}
\end{definition}

\begin{definition}\label{def:GR}
We denote by $\GR(\theta,\U\myp)$ the set of all Gibbs measures
$\g\in \mathscr{G}(\theta,\U\myp)$ such that all its correlation
functions $\kappa_{\myn\g}^{n}$ are well defined and satisfy the
\textit{Ruelle bound}, that is, for some constant $R\in \RR_+$ and
all $n\in\NN$,
\begin{equation}\label{eq:RB}
|\kappa_{\myn\g}^{n}(x_1,\dots,x_n)|\le R^{n},\qquad
(x_1,\dots,x_n)\in X^n.
\end{equation}
\end{definition}

\begin{proposition}\label{pr:k=k}
Let\/ $\g_{1},\g_{2}\in \GR(\theta,\U\myp)$ and\/
$\kappa_{\myn\g_{1}}^{n}=\kappa_{\myn\g_{2}}^{n}$ for all\/
$n\in\NN$. Then\/ $\g_{1}=\g_{2}$.
\end{proposition}
\proof
For any measure $\g\in \GR(\theta,\U\myp)$, its Laplace transform
$L_{\g}(f)$ on functions $f\in C_0(X)$ may be represented in the
form
\begin{align}
\notag L_{\g}(f)&=\int_{\varGamma _{X}}\prod_{x_i\in\gamma}
\bigl(1+(\re^{-f(x_i)}-1)\bigr)\,\g(\rd\gamma)\\
\label{eq:l1} &=1+ \int_{\varGamma _{X}}\sum_{n=1}^\infty
\sum_{\{x_{1}\myn,\dots,\myp x_{n}\} \subset \gamma}
\prod_{i=1}^n \bigl(\re^{-f(x_{i})}-1\bigr)\,\g(\rd\gamma)\\
\label{eq:l2} &=1+ \sum_{n=1}^\infty \int_{\varGamma_{X}}
\sum_{\{x_{1}\myn,\dots,\myp x_{n}\} \subset \gamma}
\prod_{i=1}^n \bigl(\re^{-f(x_{i})}-1\bigr)\,\g(\rd\gamma)\\
\label{eq:laplace} &= 1+\sum_{n=1}^{\infty}
\frac{1}{n!}\int_{X^{n}}\prod_{i=1}^n \bigl(\re^{-f(x_{i})}-1\bigr)
\,\kappa_{\myn\g}^{n}(x_{1},\dots,x_{n})\,\theta(\rd{x}_{1})\cdots\theta(\rd{x}_{n}),
\end{align}
where (\ref{eq:laplace}) is obtained from (\ref{eq:l2}) using
formula (\ref{corr-funct}). Interchanging the order of integration
and summation in (\ref{eq:l1}) is justified by the dominated
convergence theorem; indeed, using that $|f(x)|\le C_f$ on
$K_f:=\supp f$ with some $C_f>0$ and recalling that the correlation
functions $\kappa^n_{\g} $ satisfy the Ruelle bound (\ref{eq:RB}),
we see that the right-hand side of (\ref{eq:laplace}) is dominated
by
$$
1+\sum_{n=1}^{\infty} \frac{1}{n!}\,(\re^{C_f}+1)^n \mypp
R^{n}\mypp\theta(K_f)^n=\exp\bigl\{R\mypp(\re^{C_f}+1)
\mypp\theta(K_f)\bigr\}<\infty.
$$
Now, formula (\ref{eq:laplace}) implies that if measures
$\g_1,\g_2\in\GR(\theta,\U\myp)$ have the same correlation
functions, then their Laplace transforms coincide with each other,
$L_{\g_{1}}(f)=L_{\g_{2}}(f)$ for any $f\in\mathrm{M}_+(X)$, hence
$\g_1=\g_2$.
\endproof

\begin{definition}\label{def:extreme}
It is well known that $\mathscr{G}(\theta,\U\myp)$ is a convex set
\cite{Preston}. We denote by $\ext \mathscr{G}(\theta,\U\myp)$ the
set of its extreme elements, that is, those measures
$\g\in\mathscr{G}(\theta,\U\myp)$ that cannot be written as
$\g=\frac12\myp(\g_1+\g_2)$ with
$\g_1,\g_2\in\mathscr{G}(\theta,\U\myp)$ and $\g_1\ne \g_2$.
\end{definition}


Using Ruelle's equation (\ref{ruelle}), it is easy to obtain the
following result (cf.\ \cite[Corollary~2.2.6]{KunaPhD}).
\begin{proposition}\label{pr:Gibbs|cond}
Let \myp$\g\in\mathscr{G}(\theta,\U\myp)$, and let
$\varLambda\in\mathcal{B}(X)$ be a compact set.
Then the restriction of the Gibbs measure
$\g\in\mathscr{G}(\theta,\U\myp)$ onto the space
$\varGamma_\varLambda$, defined by
\begin{equation*}
\g_\varLambda(A):=\g(A\cap\varGamma_\varLambda),\qquad
A\in\mathcal{B}(\varGamma_X),
\end{equation*}
is absolutely continuous with respect to the Lebesgue--Poisson
measure $\lambda_\theta$, with the Radon--Nikodym density
$S_\varLambda:=\rd\g_\varLambda/\rd\lambda_\theta\in
L^1(\varGamma_\varLambda,\lambda_\theta)$ given by
\begin{equation}\label{eq:Psi}
S_\varLambda(\gamma)=\re^{-E(\gamma)}\int_{\varGamma_{X\setminus
\varLambda}}\re^{-E(\gamma,\myp\gamma')}\,\g(\rd\gamma'),\qquad
\gamma\in\varGamma_\varLambda.
\end{equation}
\end{proposition}

\end{document}